\documentclass[final,1p,times]{elsarticle}

\usepackage{graphicx}
\usepackage{color}
\newtheorem{thm}{Theorem}

\usepackage{amssymb,amsmath}

\usepackage{multirow}
\usepackage{stackrel}

\DeclareMathOperator{\e}{e}
\newtheorem{proposition}{Proposition}
\newdefinition{rmk}{Remark}
\newtheorem{example}{Example}
\journal{XXX}
\usepackage[all]{xy}

\begin{document}

\begin{frontmatter}



\title{Analysis of  random non-autonomous logistic-type differential equations via the Karhunen-Lo\`{e}ve expansion and the Random Variable Transformation technique}


\author[IMM]{J.-C. Cort\'es\corref{cor1}}
\ead{jccortes@imm.upv.es}

\author[deusto1,deusto2]{A. Navarro-Quiles}
\ead{annaqui@doctor.upv.es}

\author[IMM]{J.-V. Romero}
\ead{jvromero@imm.upv.es}

\author[IMM]{M.-D. Rosell\'{o}}
\ead{drosello@imm.upv.es}

\address[IMM]{Instituto Universitario de Matem\'{a}tica Multidisciplinar,\\
Universitat Polit\`{e}cnica de Val\`{e}ncia,\\
Camino de Vera s/n, 46022, Valencia, Spain}

\address[deusto1]{DeustoTech, University of Deusto, 48007 Bilbao,\\
 Basque Country, Spain}
\address[deusto2]{ Facultad de Ingenieria, Universidad de Deusto, \\
Avda.Universidades, 24, 48007, Bilbao, Basque Country, Spain.}

\cortext[cor1]{Corresponding author. Phone number: +34--963877000 (ext.~88289)}

\begin{abstract}

This paper deals with the study, from a probabilistic point of view, of logistic-type differential equations with uncertainties. We assume that  the initial condition is a random variable and the diffusion coefficient is a stochastic process. The main objective is  to obtain the first probability density function,  $f_1(p,t)$, of the solution stochastic process, $P(t,\omega)$. To achieve this goal, first  the diffusion coefficient is represented via a truncation of order $N$ of the  Karhunen-Lo\`{e}ve expansion, and second,  the  Random Variable Transformation technique is applied. In this manner, approximations, say $f_1^N(p,t)$, of  $f_1(p,t)$ are constructed. Afterwards, we rigorously prove that $f_1^N(p,t) \longrightarrow f_1(p,t)$ as $N\to \infty$ under mild conditions assumed on input data (initial condition and diffusion coefficient). Finally, three illustrative examples are shown.
\end{abstract}


\end{frontmatter}



\section{Motivation and Preliminaries}\label{section1}
The prominent role of the logistic differential equation to model problems in different settings as Biology (the dynamics of a population), Economics (the diffusion of a new technology or the growth of an economy), Engineering (the variation of physical properties subject to industrial processes), etc., has been extensively discussed and exhibited in numerous contributions (see for instance \cite{biolo,biolo2,biolo3}, \cite{econo2,logistico_tecnologia} and \cite{inge,inge2}, respectively). The logistic differential equation was first proposed by Pierre-Fran\c{c}ois Verhulst,
in his celebrated papers \cite{ver,ver2}, to overcome the shortcomings of  Malthusian's model to  study the  population growth. The main feature of Verhulst's model versus Malthus's one is the inclusion of a carrying capacity of the environment, say $M$,  which restricts the total number of individuals because resources constrains. Assuming, without loss of generality that $M=1$,  the classical logistic model is formulated via the following initial value problem (IVP)
\[
\left.\begin{array}{lcl}
p'(t)&=& a(1-p(t))p(t),\\
p(t_0)&=&p_0,
\end{array}\right\}\quad t>t_0, \quad 0<p_0<1,
\]
where $p_0$ and $p(t)$ denote the proportions of individuals at the time instants $t_0$ and $t>t_0$, respectively. This model has been thoroughly studied from different perspectives and using a number of mathematical techniques (see \cite{Islam,Pamuk}, for example).  For a fixed initial condition $p_0\in ]0,1[$,  the parameter $a>0$ stands for the reproductive parameter. This parameter depends upon complex variables including   environmental factors (weather, food, etc.),  genetic factors (birth and death rates, health, etc.), age and other influence factors whose nature is clearly random. Furthermore, the initial condition $p_0$ is often calculated via sampling techniques, thus involving randomness,  because is not feasible to  quantify   its value in an exact manner. Hence, it is more realistic to consider that  $p_0$ is a random variable (RV) rather than  a deterministic value. These reasons have motivated the study of  logistic-type differential equations with uncertainties  both in the initial condition, $p_0$, and in the reproductive parameter, $a$. Research on the logistic differential equation with randomness has  been  conducted using mainly two approaches. 

In the first one, uncertainty is introduced via stochastic processes (SPs) whose sample behaviour is very irregular (e.g., nowhere differentiability). This leads to the so-called Stochastic Differential Equations (SDEs).  For example, if stochastic perturbations (or noise) are considered by means of a Wiener process like the Brownian motion, then the rigorous treatment of the corresponding  SDE requires a special stochastic calculus whose  cornerstone result is the It\^{o} lemma \cite{Oksendall-LIBRO, Kloeden-LIBRO, Allen-LIBRO}. SDEs are formally written via stochastic differentials but rigorously analysed using Riemann-Stieltjes and It\^{o} type stochastic integrals. In this  class of SDEs, input noise is limited to   Gaussian pattern.  Some interesting contributions addressing different formulations of the logistic model or its generalizations, based upon SDEs, include \cite{Braumann_CAMWA_2008, Braumann_Fisheries_2017, Meng_CAMWA_2012,Meng_AML_2013}. 

The second approach consists of  direct randomization of input parameters (initial/boundary conditions, forcing terms and/or coefficients) by assigning them suitable probability distributions. This allows  to introduce a wider class of stochastic patterns,  including the Gaussian one, to describe uncertainties. This leads to the area of Random Differential Equations (RDEs). The so-called Random Mean Square Calculus provides a powerful tool to rigorously tackle RDEs \cite{Soong, Neckel_2013}. The study of the logistic RDE, using the Mean Square Calculus and its generalizations, can be found for instance in  \cite{MATCOM_2009, JCAM_2013_Laura}.

Additional approaches based upon SDEs/RDEs formulations  to deal with the logistic differential with uncertainty are the moment closure technique \cite{Nasell_JTP2003} and fuzzy variables \cite{Dorini_FSS2017}.

In all these contributions dealing with the logistic SDE/RDE, apart from obtaining the solution SP, say $P(t,\omega)$, a major goal is to determine its main statistical functions, namely, the mean function, $\mathbb{E}\left[P(t,\omega)\right]$, and the variance function, $\mathbb{V}\left[P(t,\omega)\right]$. However, a more ambitious target is the computation of its first probability density  function (1-PDF), $f_1(p,t)$, since via its integration one can compute  all the one-dimensional statistical moment functions, 
\begin{equation}\label{moments}
\mathbb{E} [\left(P(t,\omega)\right)^k ]=\int_{-\infty}^{\infty} p^k f_1(p,t)\,\mathrm{d}p,\quad k=1,2,\ldots, 
\end{equation}
and, in particular, the mean, $\mathbb{E}\left[P(t,\omega)\right]$, and the variance,
\begin{eqnarray}\label{varianza_formula}
\mathbb{V}\left[P(t,\omega)\right] =\mathbb{E}[(P(t,\omega))^2]-\mathbb{E}^2 [P(t,\omega)],
\end{eqnarray}
as well as the probability that the population lies in a set of specific interest, say $[p_1,p_2]$,
\[
\mathbb{P}\left[p_1 \leq P(t,\omega) \leq p_2\right]=\int_{p_1}^{p_2} f_1(p,t)\, \mathrm{d}p.
\]
On the one hand, and to the best of our knowledge,  the computation of the 1-PDF of the solution SP of the logistic RDE has been addressed in the two following recent contributions \cite{Dorini_CNSNS2016, Dorini_CAM2016}, but only in the case that coefficients do not depend on time, i.e., for the so-called  random autonomous logistic differential equation.  On the other hand, the authors have recently  obtained the 1-PDF of the solution SP to the random non-autonomous first-order linear homogeneous differential equation by combing the Karhunen-Lo\`{e}ve expansion and the  Random Variable Transformation technique \cite{KL1}. Aimed by this latter result,  the goal of this paper is to extend the analysis performed in \cite{KL1} to the random logistic differential equation assuming that the reproductive parameter is a time-dependent SP, instead of being a RV, and further assuming that the initial condition is a RV. In this manner, here we will deal with the general case from a probabilistic standpoint.

Specifically, we will consider the following random IVP  
\begin{equation}\label{logistic_problem_random}
\left.\begin{array}{lcl}
P'(t,\omega)&=&A(t,\omega)(1-P(t,\omega))P(t,\omega),\\
P(t_0,\omega)&=&P_0(\omega),
\end{array}\right\}\quad t>t_0.
\end{equation}
where $A(t,\omega)$ is a SP, $P_0(\omega)$ is a bounded absolutely  continuous  RV   $ P_0: \Omega \longrightarrow [p_{0,1}, p_{0,2}]\subset ]0,1[$ both satisfying certain hypotheses that will be specified later. These random inputs are defined on a common complete probability space $(\mathbb{P},\mathcal{F},\Omega)$. 

Although in the biological setting, the reproductive coefficient in the  logistic differential equation is naturally positive, for the sake of generality,  in the subsequent analysis we will deal with the general case where  $A(t,\omega)$ is not necessarily positive. In this manner, our study will be also valid in more general terms.

Therefore, the main goal of this paper is to compute the 1-PDF of the solution SP of the random IVP \eqref{logistic_problem_random}. To reach this objective, we will combine two important results, namely,   the Karhunen-Lo\`eve expansion (KLE) and the  Random Variable Transformation (RVT) method. The former result will be applied to represent the  coefficient $A(t,\omega)$ via an expansion of a denumerable set of zero-mean, unit variance and uncorrelated  RVs. By truncating, up to certain order $N$ the KLE of $A(t,\omega)$,  we  will then apply the RVT method  to determine approximations, $f_1^N(p,t)$ to the exact 1-PDF $f_1(p,t)$, for each $(p,t)$ fixed.  Afterwards, we will prove the convergence of $f_1^N(p,t) \longrightarrow f_1(p,t)$ as $N\to \infty$ assuming certain hypotheses on input data $P_0(\omega)$ and $A(t,\omega)$ that will be specified later on. For the sake of completeness, down below we state both the RVT technique as the KLE.

\begin{thm}[Random Variable Transformation technique]\label{RVT}\cite{Soong}
Let  $\mathbf{X}(\omega)=\left( X_1(\omega),\ldots,X_m(\omega) \right)^{\mathsf{T}}$ and $\mathbf{Y}(\omega)=\left( Y_1(\omega),\ldots,Y_m(\omega)\right) ^{\mathsf{T}}$  be two $m$-dimensional absolutely continuous random vectors defined on a complete probability space $(\Omega,\mathfrak{F},\mathbb{P})$. Let $\mathbf{r}: \mathbb{R}^m \rightarrow \mathbb{R}^m$ be a one-to-one deterministic transformation of $\mathbf{X}(\omega)$ into $\mathbf{Y}(\omega)$, i.e., $\mathbf{Y}(\omega)=\mathbf{r}(\mathbf{X}(\omega))$, $\omega\in \Omega$. Assume that $\mathbf{r}$ is a continuous mapping and has continuous partial derivatives with respect to each component $x_i$, $1\leq i \leq m$. Then, if $f_{\mathbf{x}}(x_1,\ldots,x_m)$ denotes the joint probability density function of the vector $\mathbf{X}(\omega)$, and $\mathbf{s}=\mathbf{r}^{-1}=(s_1(y_1,\ldots,y_m),\ldots,s_m(y_1,\ldots,y_m))$ represents the inverse mapping of $\mathbf{r}=(r_1(x_1,\ldots,x_m),\ldots,r_m(x_1,\ldots,x_m))$, the joint probability density function of the random vector $\mathbf{Y}(\omega)$ is given by
\[
f_{\mathbf{Y}}(y_1,\ldots,y_m)=f_{\mathbf{X}}\left(s_1(y_1,\ldots,y_m),\ldots,s_m(y_1,\ldots,y_m)\right) \left| \mathcal{J}_m\right|,
\]
where $\left| \mathcal{J}_m\right|$, which is assumed to be different from zero, denotes the absolute value of the Jacobian defined by the following determinant
\[
\mathcal{J}_m
=
\det
\left[
\begin{array}{ccc}
\displaystyle \frac{\partial s_1(y_1,\ldots, y_m)}{\partial y_1} & \cdots & \displaystyle \frac{\partial s_m(y_1,\ldots, y_m)}{\partial y_1}\\
\vdots & \ddots & \vdots\\
\displaystyle \frac{\partial s_1(y_1,\ldots, y_m)}{\partial y_m} & \cdots & \displaystyle \frac{\partial s_m(y_1,\ldots, y_m)}{\partial y_m}\\
\end{array}
\right]
\,.
\] 
\end{thm}

It is worthy highlighting that in the context of RDEs, the RVT technique has been successfully applied to determine  the 1-PDF of the solution stochastic process of relevant problems, formulated via RDEs, that appear in different areas  \cite{Hussein-Selim2012,Slama-Selim2017,Dorini2011_MATCOM,Santos-Dorini-Cunha-AMC-2010}.

\begin{thm}[$\mathrm{L}^2$ convergence of Karhunen-Lo\`eve]\label{KLE}\cite[p.~202]{Powell-LIBRO}
Consider a mean square integrable continuous time stochastic process  $X\equiv \{X(t,\omega): t\in \mathcal{T},\omega \in \Omega\}$, i.e.,  $X\in \mathrm{L}^2(\Omega , \mathrm{L}^2(\mathcal{T}))$ being $\mu_X (t)$ and $c_X(s,t)$ its  mean and covariance functions, respectively. Then, 
\begin{equation}\label{KLEX}
X(t,\omega)=\mu_X(t)+\sum_{j=1}^{\infty} \sqrt{\nu_j} \, \phi_j(t)\,\xi_j(\omega),\quad \omega \in \Omega,
\end{equation}
where, this sum converges in $\mathrm{L}^2(\Omega, \mathrm{L}^2(\mathcal{T}))$,
\[
\xi_j(\omega):=\frac{1}{\sqrt{\nu_j}}\left< X(t,\omega)-\mu_X(t),\phi_j(t)\right>_{\mathrm{L}^2(\mathcal{T})},
\]
$\{(\nu_j,\phi_j(t)): j\geq 1\}$ denote, respectively, the eigenvalues with $\nu_1 \geq \nu_2 \geq \cdots \geq 0$ and  eigenfunctions of the following integral    operator $\mathfrak{C}$ 
\[
(\mathfrak{C}f)(t):=\int_{\mathcal{T}} c_X(s,t)f(s)\, \mathrm{d} s, \quad f\in \mathrm{L}^2(\mathcal{T}),
\]
associated to the covariance function $c_X(s,t)$. RVs $\xi_j (\omega)$ have zero mean ($\mathbb{E} [\xi_j(\omega)]=0$), unit variance ($\mathbb{V} [\xi_j(\omega)]=1$) and are pairwise uncorrelated ($\mathbb{E}  [\xi_j(\omega) \xi_k(\omega)]=\delta_{jk}$). Furthermore, if $X(t,\omega)$ is Gaussian, then $\xi_j(\omega) \sim \mathrm{N}(0,1)$ are independent and identically distributed.
\end{thm}

The space $(\mathrm{L}^2(\Omega , \mathrm{L}^2(\mathcal{T})), \left\| \cdot \right\|_{\mathrm{L}^2(\Omega , \mathrm{L}^2(\mathcal{T}))})$ introduced in Th.~ \ref{KLE} corresponds to the set of square integrable SPs, $X(t,\omega)$, defined in a set $\mathcal{T}\subset \mathbb{R}$, i.e., $\int_{\mathcal{T}} \mathbb{E} \left[ |X(t,\omega)|^2 \right]  \mathrm{d} t < \infty$ (see \cite{Powell-LIBRO}) with the norm
\begin{equation}\label{norma_L_2}
 \left\| X(t,\omega) \right\|_{\mathrm{L}^2(\Omega , \mathrm{L}^2(\mathcal{T}))}
=
\left(
\int_{\mathcal{T}} \mathbb{E} \left[ |X(t,\omega)|^2 \right]  \mathrm{d} t 
\right)^{1/2}<\infty.
\end{equation}

We finish this section by stating two results that will be used in the last example of Section \ref{sec_ejemplos}.
\begin{proposition}(\cite[Th.~8, p.~92]{Stirzaker}) \label{propo_independencia}
Let $\{\xi_i(\omega): 1\leq i \leq N\}$ be independent real random variables defined in a common probability space $(\Omega,\mathcal{F},\mathbb{P})$. Let $f_i: \Omega \longrightarrow \mathbb{R}$, $1\leq i \leq N$, be Borel measurable functions. Then, $\{f_i(\xi_i(\omega)): 1\leq i \leq N\}$ are independent random variables.
\end{proposition}
\begin{proposition}(\cite[p.~21]{Massart}) \label{desig_clave}
Let $\xi(\omega)$ be a random variable defined in a probability space $(\Omega,\mathcal{F},\mathbb{P})$ such that $\mathbb{E}[\xi(\omega)]=0$ and $\mathbb{P}[\{  \omega: \alpha \leq \xi(\omega) \leq \beta]=1$. Then,
\[
\mathbb{E}\left[ \e^{ \lambda \xi(\omega)}\right]
\leq
\e^{\frac{\lambda^2 (\beta - \alpha)^2}{8}},\quad \lambda \in \mathbb{R}. 
\]
\end{proposition}
For notational convenience, throughout this paper the exponential function will be written  by $\e^x$ or $\exp(x)$, interchangeably.

\section{Main result: Computing the 1-PDF of the solution stochastic process}
This section is firstly addressed to construct approximations, $f_1^N(p,t)$, of the 1-PDF, $f_1(p,t)$, of the solution SP, $P(t,\omega)$, to the random non-autonomous  IVP  \eqref{logistic_problem_random}, and secondly,  to prove that these approximations are convergent, i.e.  $f_1^N(p,t) \longrightarrow f_1(p,t)$ as $N \to \infty$,  assuming mild conditions  on the random inputs  $P_0(\omega)$  and $A(t,\omega)$. 

As the construction of the aforementioned approximation and the proof of its convergence follows a rather technical process, for the sake of clarity we first present an overview of the main ideas that  that will be applied to achieve  these two goals.

Regarding the construction of the approximation, first we will consider the formal randomization of solution, $P(t,\omega)$, of the logistic model (see below expression \eqref{logistic_solution}), which depends on the stochastic process $A(t)$. Second, we will represent $P(t,\omega)$ in terms of the KLE of  $A(t,\omega)$ (see expression \eqref{logistic_solution_random}). Next, we will truncate this KLE, say $A_N(t)$, so that we will obtain a formal approximation, $P_N(t,\omega)$, of the solution stochastic process $P(t,\omega)$ (see expression  \eqref{P_N}). Then we will apply the RVT method, stated in Th.~\ref{RVT}, to obtain  the  1-PDF, $f_1^N(p,t)$,  of  $P_N(t,\omega)$ (see expression \eqref{1PDF}). To legitimate this approach, we will assume two hypotheses, that will be denoted by \textbf{H1} and \textbf{H2}. As we will comment later on, \textbf{H1} allows us to assure that $P(t,\omega)$ is well defined from a probabilistic standpoint, while \textbf{H2} guarantees the initial condition, $P_0(\omega)$, and the random variables involved in the KLE (see Th.~\ref{KLE}), say $(\xi_1(\omega),\ldots,\xi_N(\omega)):=\boldsymbol{\xi}_{N}(\omega)$, possess respective PDFs. Additionally, all these RVs are assumed to be independent which is a natural assumption in our stochastic setting.

To proof that $f_1^N(p,t) \longrightarrow f_1(p,t)$ as $N \to \infty$, we will apply the classical Cauchy condition assuming two further hypotheses \textbf{H3} and \textbf{H4}. Hypothesis \textbf{H3} assumes that the PDF of the initial condition, $f_{P_0}(p_0)$, is Lipschitz. As usual, this assumption permits transferring the behaviour of the increment in the range of $f_{P_0}$, that appears when the Cauchy condition is applied, in terms of the information of its domain which involves random variables $P_0(\omega)$ and $\xi_i(\omega)$, $1\leq i \leq N$. This strategy allows us to take advantage of hypothesis \textbf{H4} which is related to the growth  of the moment-generating function of the KLE of the integral stochastic process $K_N(t,\boldsymbol{\xi}_{N}(\omega))\int_{t_0}^t A_N(s,\omega)\, \mathrm{d}s$ (see \eqref{K_N}).

After following the approach previously described, we will summarize our conclusions in a theorem (see  Theorem \ref{main_result}).

Hereinafter, we will assume that  $t\in \mathcal{T}=[t_0,T]$, $T>t_0$. Motivated by its deterministic counterpart, it is easy to check that a formal solution SP of random IVP \eqref{logistic_problem_random} is given by
\begin{equation}\label{logistic_solution}
P(t,\omega)=\frac{1}{1+\exp{ \left( \int_{t_0}^t -A(s,\omega)\,\mathrm{d}s \right)}\left(-1+\frac{1}{P_0(\omega)}\right)}, \quad t\in \mathcal{T}, \,\, \omega \in \Omega.
\end{equation}
We will assume that both random inputs satisfy the following hypothesis:
\[
\textbf{H1}:
P_0: \Omega \longrightarrow [p_{0,1}, p_{0,2}] \subset ]0,1[ \,\, \text{and}\,\,\, A(t,\omega)\in \mathrm{L}^2(\Omega,\mathrm{L}^2(\mathcal{T})).
\]
\begin{rmk}\label{P_N_entre_0y1}
Notice that the first part of this hypothesis guarantees that the denominator of $P(t,\omega)$ is nonzero with probability 1, thus $P(t,\omega)$ is  well-defined almost everywhere (a.e.) regardless the sign of $A(t,\omega)$. Moreover, from \eqref{logistic_solution} it is clear that $P(t,\omega)\neq 0$ a.e. Thus, $P(t,\omega)\in ]0,1[$ for all $t\in [t_0,T]$ and $\omega \in \Omega$ a.e. 
\end{rmk} 
On the one hand, notice that as we are assuming  the initial condition is bounded, then it is a second-order RV, i.e. $\mathbb{E} [(P_0(\omega))^2]<+\infty$ (hence having finite variance) and as the time interval $\mathcal{T}$ is closed and bounded, then $P_0(\omega)\in \mathrm{L}^2(\Omega,\mathrm{L}^2(\mathcal{T}))$. On the other hand, since $A(t,\omega) \in \mathrm{L}^2(\Omega,\mathrm{L}^2(\mathcal{T}))$,  it can be represented via the KLE given in \eqref{KLEX}. Then, the formal solution SP, given in \eqref{logistic_problem_random}, can be written as 
\begin{equation}\label{logistic_solution_random}
P(t,\omega)=\displaystyle \frac{1}{1+\exp \left(-\displaystyle \int_{t_0}^t \left(  \mu_A(t)+\sum_{j=1}^{\infty} \sqrt{\nu_j} \, \phi_j(t)\,\xi_j(\omega) \right)\, \mathrm{d}s\right)\left(-1+\dfrac{1}{P_0(\omega)}\right)}, \quad t\in \mathcal{T}, \,\, \omega \in \Omega.
\end{equation}

Now, we will use this fact together with the RVT technique to construct the approximations $f_1^N(p,t)$. With this aim, let us consider the truncation of  order $N$ of the KLE to  $A(t,\omega)$
\[
A_N(t,\omega)=\mu_A(t)+\sum_{j=1}^{N} \sqrt{\nu_j}\, \phi_j(t)\,\xi_j(\omega).
\]
Thus, according to \eqref{logistic_solution_random},  the following formal approximation of the truncated solution SP is obtained
\begin{equation}\label{P_N}
P_N(t,\omega)=\frac{1}{1+\exp \left( \displaystyle-\int_{t_0}^t \left(  \mu_A(s)+\sum_{j=1}^{N} \sqrt{\nu_j}\, \phi_j(s)\,\xi_j(\omega) \right) \mathrm{d}s\right)\left(\displaystyle -1+\frac{1}{P_0(\omega)}\right)}.
\end{equation}
Now, in order to apply the RVT technique, we will assume the following hypothesis
\[
\textbf{H2}:
\begin{array}{cc}
P_0(\omega), \xi_i(\omega),\,\, 1\leq i \leq N, \,\, \text{are absolutely continuous RVs.}\\
P_0(\omega),\,\,\, \boldsymbol{\xi}_{N}(\omega)=(\xi_1(\omega),\dots, \xi_N(\omega)) \,\, \text{are independent}\\
\text{with PDFs} \,\, 
f_{P_0}(p_0)\,\,\text{and} \,\, 
f_{\boldsymbol{\xi}_{N}}(\xi_1,\dots, \xi_N),\,\, \text{respectively}.
\end{array}
\]
Then,  we define the following one-to-one transformation $\mathbf{r}: \mathbb{R}^{N+1} \longrightarrow \mathbb{R}^{N+1}$ componentwise
\[
\begin{array}{ccccl}
y_1&=&r_1(p_0,\xi_1,\dots,\xi_N)&=& \left(1+\exp\left(\displaystyle-\int_{t_0}^t \left(  \mu_A(s)+\sum_{j=1}^{N} \sqrt{\nu_j}\, \phi_j(s)\,\xi_j \right)  \mathrm{d}s\right)\left(\displaystyle -1+\frac{1}{p_0}\right)\right)^{-1},\\
\\
y_2 &= &r_2(p_0,\xi_1,\dots,\xi_N)&=& \xi_1,\\
\vdots & & \vdots\\
y_{N+1} &= &r_{N+1}(p_0,\xi_1,\dots,\xi_N)&=& \xi_N.
\end{array}
\]
The inverse of mapping $\mathbf{r}$ is   $\mathbf{s}: \mathbb{R}^{N+1} \longrightarrow \mathbb{R}^{N+1}$,  whose components are
\[
\begin{array}{ccccl}
p_0&=&s_1(y_1,\dots, y_{N+1})&=& \displaystyle \frac{y_1 \exp\left(\displaystyle-\int_{t_0}^t \left(  \mu_A(s)+\sum_{j=1}^{N} \sqrt{\nu_j}\, \phi_j(s)\,y_{j+1} \right) \mathrm{d}s\right)}{1+y_1\left(-1+\exp\left(\displaystyle-\int_{t_0}^t \left(  \mu_A(s)+\sum_{j=1}^{N} \sqrt{\nu_j}\, \phi_j(s)\,y_{j+1} \right) \mathrm{d}s\right)\right)},\\
\\
\xi_1 &=&s_2(y_1,\dots, y_{N+1})&=& y_2,\\
\vdots & & \vdots & & \vdots\\
\xi_N &=&s_{N+1}(y_1,\dots, y_{N+1})&=& y_{N+1}.
\end{array}
\]
Moreover, the absolute value of the Jacobian of the inverse mapping $\mathbf{s}$ is nonzero since
\[
\left|J \right|=\left| \frac{\partial s_1}{\partial y_1} \right|=\frac{\exp\left(\displaystyle-\int_{t_0}^t \left(  \mu_A(s)+\sum_{j=1}^{N} \sqrt{\nu_j}\, \phi_j(s)\,y_{j+1} \right) \mathrm{d}s\right)}{\left(1+y_1\left(-1+\exp\left(\displaystyle-\int_{t_0}^t \left(  \mu_A(s)+\sum_{j=1}^{N} \sqrt{\nu_j}\, \phi_j(s)\,y_{j+1} \right) \mathrm{d}s\right)\right)\right)^2}\neq 0.
\]
Therefore, applying  Th.~\ref{RVT}, the PDF of the random vector $(Y_1(\omega),\dots, Y_{N+1}(\omega))$ is given by
\[
\begin{array}{lcl}
f_{Y_1,\dots,Y_{N+1}} (y_1,\dots, y_{N+1})&=& \displaystyle f_{P_0}\left( \frac{y_1 \exp\left(\displaystyle-\int_{t_0}^t \left(  \mu_A(s)+\sum_{j=1}^{N} \sqrt{\nu_j}\, \phi_j(s)\,y_{j+1} \right) \mathrm{d}s\right)}{1+y_1\left(-1+\exp\left(\displaystyle-\int_{t_0}^t \left(  \mu_A(s)+\sum_{j=1}^{N} \sqrt{\nu_j}\, \phi_j(s)\,y_{j+1} \right) \mathrm{d}s\right)\right)} \right)
\\
\\
& & \times
f_{\boldsymbol{\xi}_N}\left( y_2,\dots, y_{N+1}\right)
\\
\\
& & \times \displaystyle \frac{\exp\left(\displaystyle-\int_{t_0}^t \left(  \mu_A(s)+\sum_{j=1}^{N} \sqrt{\nu_j}\, \phi_j(s)\,y_{j+1} \right) \mathrm{d}s\right)}{\left(1+y_1\left(-1+\exp\left(\displaystyle-\int_{t_0}^t \left(  \mu_A(s)+\sum_{j=1}^{N} \sqrt{\nu_j}\, \phi_j(s)\,y_{j+1} \right) \mathrm{d}s\right)\right)\right)^2}.
\end{array}
\]
Finally, marginalizing with respect to the random vector $(Y_2(\omega),\dots, Y_{N+1}(\omega))=\boldsymbol{\xi}_N (\omega)$  and ta\-king $t\geq t_0$ arbitrary, we obtain the following explicit expression of the 1-PDF of the truncated solution SP, $P_N(t,\omega)$,
\begin{equation}\label{1PDF}
\begin{array}{lcl}
f_1^N(p,t)&=& \displaystyle \int_{\mathbb{R}^N} f_{P_0}\left( \frac{p \exp\left(\displaystyle-\int_{t_0}^t \left(  \mu_A(s)+\sum_{j=1}^{N} \sqrt{\nu_j}\, \phi_j(s)\,\xi_{j} \right) \mathrm{d}s\right)}{1+p\left(-1+\exp\left(\displaystyle-\int_{t_0}^t \left[  \mu_A(s)+\sum_{j=1}^{N} \sqrt{\nu_j}\, \phi_j(s)\,\xi_{j} \right] \mathrm{d}s\right)\right)}\right)\\
\\
& & \times
\displaystyle   f_{\boldsymbol{\xi}_N}\left( \xi_1,\dots, \xi_N\right)
\\
\\
& & \times  \displaystyle \frac{\exp\left(\displaystyle-\int_{t_0}^t \left(  \mu_A(s)+\sum_{j=1}^{N} \sqrt{\nu_j}\, \phi_j(s)\,\xi_{j} \right) \mathrm{d}s\right)}{\left(1+p\left(-1+\exp\left(\displaystyle-\int_{t_0}^t \left(  \mu_A(s)+\sum_{j=1}^{N} \sqrt{\nu_j}\, \phi_j(s)\,\xi_{j} \right) \mathrm{d}s\right)\right)\right)^2} \mathrm{d}\xi_N \cdots \mathrm{d}\xi_1.
\end{array}
\end{equation}

In our subsequent analysis,   we  will provide  conditions in order to guarantee the following convergence
\[
\lim_{N\to +\infty} f_1^N(p,t)=f_1(p,t), \quad  \forall (p,t)\in \mathcal{J} \times  [t_0,T]\, \text{fixed, being} \,\, \mathcal{J} \subset ]0,1[ \,\, \text{bounded}.
\]
Hereinafter, $\hat{p}$ will denote a lower bound of $\mathcal{J}$, i.e., $\hat{p} \in ]0,1[$ such that $\hat{p}<p$ for all $p\in \mathcal{J}$. We will prove this convergence by applying the Cauchy condition, i.e., for $\epsilon>0$ fixed, there exist $n_0$ (independent of $(p,t)$), such as 
\[
\left| f_1^N(p,t)-f_1^M(p,t)\right| < \epsilon, \quad \forall (p,t)\in \mathcal{J} \times  [t_0,T]\, \text{fixed, being} \,\, \mathcal{J} \subset ]0,1[ \,\, \text{bounded},\,\, \forall N,M\geq n_0.
\]

In order to simplify notation, henceforth we will denote
\begin{equation}\label{K_N}
K_N(t,\boldsymbol{\xi}_{N}(\omega))=\int_{t_0}^{t} \left( \mu_A(s)+\sum_{j=1}^{N} \sqrt{\nu_j}\, \phi_j(s) \, \xi_{j}(\omega) \right) \mathrm{d}s.
\end{equation}
Thus, expression \eqref{1PDF} can be rewritten as
\begin{equation}\label{1PDF_K_N}
f_1^N(p,t)= \displaystyle \int_{\mathbb{R}^N} f_{P_0}\left( \frac{p \e^{-K_N(t,\boldsymbol{\xi}_{N})}}{1+p\left(-1+\e^{-K_N(t,\boldsymbol{\xi}_{N})}\right)}\right)\frac{\e^{-K_N(t,\boldsymbol{\xi}_{N})}}{\left(1+p\left(-1+\e^{-K_N(t,\boldsymbol{\xi}_{N})}\right)\right)^2} f_{\boldsymbol{\xi}_N}(\xi_1,\dots, \xi_N) \mathrm{d}\xi_N\cdots \mathrm{d}\xi_1.
\end{equation}
Throughout the proof the Steps (I)-(III) will be done. For the sake of clarify, we legitimate these steps later on.   Let $(p,t)\in \mathcal{J} \times  [t_0,T]$ fixed, being $\mathcal{J} \subset ]0,1[$  bounded, and without loss of generality, take $N>M$. Then, 

\[
\begin{array}{l}
\left| f_1^N(p,t)-f_1^M(p,t) \right|\\
\\
=\left| \displaystyle \int_{\mathbb{R}^N} f_{P_0}\left( \frac{p \e^{-K_N(t,\boldsymbol{\xi}_{N})}}{1+p\left(-1+\e^{-K_N(t,\boldsymbol{\xi}_{N})}\right)}\right) \frac{\e^{-K_N(t,\boldsymbol{\xi}_{N})}}{\left(1+p\left(-1+\e^{-K_N(t,\boldsymbol{\xi}_{N})}\right)\right)^2}f_{\boldsymbol{\xi}_N}(\xi_1,\dots, \xi_N)\mathrm{d}\xi_N\cdots \mathrm{d}\xi_1\right.\\
\\
\quad \left. \displaystyle -\int_{\mathbb{R}^M} f_{P_0}\left( \frac{p \e^{-K_M(t,\boldsymbol{\xi}_{M})}}{1+p\left(-1+\e^{-K_M(t,\boldsymbol{\xi}_{M})}\right)}\right) \frac{\e^{-K_M(t,\boldsymbol{\xi}_{M})}}{\left(1+p\left(-1+\e^{-K_M(t,\boldsymbol{\xi}_{M})}\right)\right)^2}f_{\boldsymbol{\xi}_M}(\xi_1,\dots, \xi_M)\mathrm{d}\xi_M\cdots \mathrm{d}\xi_1\right|\\
\\
\displaystyle   \stackrel{(\mathrm{I})}{=} \left| \int_{\mathbb{R}^N} \left[ f_{P_0}\left( \frac{p \e^{-K_N(t,\boldsymbol{\xi}_{N})}}{1+p\left(-1+\e^{-K_N(t,\boldsymbol{\xi}_{N})}\right)}\right)\frac{\e^{-K_N(t,\boldsymbol{\xi}_{N})}}{\left(1+p\left(-1+\e^{-K_N(t,\boldsymbol{\xi}_{N})}\right)\right)^2} \right. \right.\\
\\
\qquad \displaystyle \left.\left. -  f_{P_0}\left( \frac{p \e^{-K_M(t,\boldsymbol{\xi}_{M})}}{1+p\left(-1+\e^{-K_M(t,\boldsymbol{\xi}_{M})}\right)}\right)\frac{\e^{-K_M(t,\boldsymbol{\xi}_{M})}}{\left(1+p\left(-1+\e^{-K_M(t,\boldsymbol{\xi}_{M})}\right)\right)^2}\right]f_{\boldsymbol{\xi}_N}(\xi_1,\dots, \xi_N)\mathrm{d}\xi_N\cdots \mathrm{d}\xi_1\right|\\
\\
\displaystyle \leq \int_{\mathbb{R}^N}   \left| \left[ f_{P_0}\left( \frac{p \e^{-K_N(t,\boldsymbol{\xi}_{N})}}{1+p\left(-1+\e^{-K_N(t,\boldsymbol{\xi}_{N})}\right)}\right)\frac{\e^{-K_N(t,\boldsymbol{\xi}_{N})}}{\left(1+p\left(-1+\e^{-K_N(t,\boldsymbol{\xi}_{N})}\right)\right)^2}\right.\right.\\
\\
\qquad \displaystyle - f_{P_0}\left( \frac{p \e^{-K_N(t,\boldsymbol{\xi}_{N})}}{1+p\left(-1+\e^{-K_N(t,\boldsymbol{\xi}_{N})}\right)}\right)\frac{\e^{-K_M(t,\boldsymbol{\xi}_{M})}}{\left(1+p\left(-1+\e^{-K_M(t,\boldsymbol{\xi}_{M})}\right)\right)^2}\\
\\
\qquad \displaystyle + f_{P_0}\left( \frac{p \e^{-K_N(t,\boldsymbol{\xi}_{N})}}{1+p\left(-1+\e^{-K_N(t,\boldsymbol{\xi}_{N})}\right)}\right)\frac{\e^{-K_M(t,\boldsymbol{\xi}_{M})}}{\left(1+p\left(-1+\e^{-K_M(t,\boldsymbol{\xi}_{M})}\right)\right)^2}
\\
\\
\qquad \displaystyle \left.\left.-  f_{P_0}\left( \frac{p \e^{-K_M(t,\boldsymbol{\xi}_{M})}}{1+p\left(-1+\e^{-K_M(t,\boldsymbol{\xi}_{M})}\right)}\right)\frac{\e^{-K_M(t,\boldsymbol{\xi}_{M})}}{\left(1+p\left(-1+\e^{-K_M(t,\boldsymbol{\xi}_{M})}\right)\right)^2}\right]\right|f_{\boldsymbol{\xi}_N}(\xi_1,\dots, \xi_N)\mathrm{d}\xi_N\cdots \mathrm{d}\xi_1\\
\\
\displaystyle \leq \int_{\mathbb{R}^N}   \left[ \underbrace{ f_{P_0}\left( \frac{p \e^{-K_N(t,\boldsymbol{\xi}_{N})}}{1+p\left(-1+\e^{-K_N(t,\boldsymbol{\xi}_{N})}\right)} \right) }_{(1)} \underbrace{\left| \frac{\e^{-K_N(t,\boldsymbol{\xi}_{N})}}{\left(1+p\left(-1+\e^{-K_N(t,\boldsymbol{\xi}_{N})}\right)\right)^2}-\frac{\e^{-K_M(t,\boldsymbol{\xi}_{M})}}{\left(1+p\left(-1+\e^{-K_M(t,\boldsymbol{\xi}_{M})}\right)\right)^2}\right| }_{(2)}\right.
\end{array}
\]
\[
\begin{array}{l}
\displaystyle +\left.  \underbrace{\left|  f_{P_0}\left( \frac{p \e^{-K_N(t,\boldsymbol{\xi}_{N})}}{1+p\left(-1+\e^{-K_N(t,\boldsymbol{\xi}_{N})}\right)}\right) -f_{P_0}\left( \frac{p \e^{-K_M(t,\boldsymbol{\xi}_{M})}}{1+p\left(-1+\e^{-K_M(t,\boldsymbol{\xi}_{M})}\right)}  \right)\right|}_{(3)} \underbrace{ \frac{\e^{-K_M(t,\boldsymbol{\xi}_{M})}}{\left(1+p\left(-1+\e^{-K_M(t,\boldsymbol{\xi}_{M})}\right)\right)^2}}_{(4)}\right] f_{\boldsymbol{\xi}_N}(\xi_1,\dots, \xi_N)\mathrm{d}\xi_N\cdots \mathrm{d}\xi_1\\
\\
\displaystyle \stackrel{(\mathrm{II})}{<}\int_{\mathbb{R}^N}\left[ \left(L_{f_{P_0}}(1+p_{0,1})+F_0\right)  C_g\vert K_N(t,\boldsymbol{\xi}_{N})-K_M(t,\boldsymbol{\xi}_{M})\vert \right.\\
\\
\displaystyle \quad   \left.+L_{f_{P_0}}C_h\vert K_N(t,\boldsymbol{\xi}_{N})-K_M(t,\boldsymbol{\xi}_{M})\vert \frac{\e^{K_M(t,\boldsymbol{\xi}_{M})}}{\hat{p}^2}\right]f_{\boldsymbol{\xi}_N}(\xi_1,\dots, \xi_N)\mathrm{d}\xi_N\dots \mathrm{d}\xi_1\\
\\
\displaystyle = \left(L_{f_{P_0}}(1+p_{0,1})+F_0 \right)  C_g\mathbb{E}\left[\vert K_N(t,\boldsymbol{\xi}_{N}(\omega))-K_M(t,\boldsymbol{\xi}_{M}(\omega))\vert \right]\\
\\
\displaystyle \quad +\frac{L_{f_{P_0}}C_h}{\hat{p}^2}\mathbb{E}\left[\vert K_N(t,\boldsymbol{\xi}_{N}(\omega))-K_M(t,\boldsymbol{\xi}_{M}(\omega))\vert \e^{K_M(t,\boldsymbol{\xi}_{M}(\omega))}\right]\\
\\
\displaystyle \leq \left(L_{f_{P_0}}(1+p_{0,1})+F_0 \right)  C_g\mathbb{E}\left[\vert K_N(t,\boldsymbol{\xi}_{N}(\omega))-K_M(t,\boldsymbol{\xi}_{M}(\omega))\vert^2 \right]^{1/2}
\\
\\
\displaystyle \quad+\frac{L_{f_{P_0}}C_h}{\hat{p}^2}\mathbb{E}\left[\vert K_N(t,\boldsymbol{\xi}_{N}(\omega))-K_M(t,\boldsymbol{\xi}_{M}(\omega))\vert^2 \right]^{1/2}\mathbb{E}\left[ \e^{2K_M(t,\boldsymbol{\xi}_{M}(\omega))}\right]^{1/2} \\
\\
\displaystyle =\left(\left(L_{f_{P_0}}(1+p_{0,1})+F_0 \right)  C_g+\frac{L_{f_{P_0}}C_h}{\hat{p}^2}\mathbb{E}\left[ \e^{2K_M(t,\boldsymbol{\xi}_{M}(\omega))}\right]^{1/2}\right)\mathbb{E}\left[\vert K_N(t,\boldsymbol{\xi}_{N}(\omega))-K_M(t,\boldsymbol{\xi}_{M}(\omega))\vert^2 \right]^{1/2} \\
\\
\stackrel{(\mathrm{III})}{\leq}\alpha\, \mathbb{E}\left[\vert K_N(t,\boldsymbol{\xi}_{N}(\omega))-K_M(t,\boldsymbol{\xi}_{M}(\omega))\vert^2 \right]^{1/2}.
\end{array}
\]

Notice that in the last steps (penultimate equality and inequality), we have applied the definition of the expectation and the Cauchy--Schwarz inequality for expectations, respectively.

Now, we justify the above Steps (I)-(III):\\
\textbf{Step \textrm{(I)}}: As $N>M$, here we have used that the  PDF $f_1^M(p,t)$ can be expressed in terms of   PDF $f_1^N(p,t)$, by marginalizing this latter function with respect to $\xi_{M+1}(\omega), \ldots,\xi_N(\omega)$, that is  by introducing the co\-rres\-pon\-ding $N-M$-fold integration. \\
\textbf{Step (II)}: The following hypothesis will be assumed to legitimate bounds of this step.
 \[
\textbf{H3}:\begin{array}{c}
f_{P_0}(p_0)\,\,
\text{is Lipschitz continuous in its domain, i.e.,}\,\,\\
\exists\, L_{f_{P_0}}>0:\, |f_{P_0}(p^*)-f_{P_0}(p^{**})| \leq L_{f_{P_0} }|p^*-p^{**}|,\quad \forall p^*,p^{**}\in [p_{0,1}, p_{0,2}]\subset]0,1[.
\end{array}
\]
Using this assumption, we get the bound corresponding to term (1). Let $F_0= f_{P_0}(p_{0,1})$, by hypothesis \textbf{H3}:
\[
\begin{array}{lcl}
\displaystyle f_{P_0}\left( \frac{p \e^{-K_N(t,\boldsymbol{\xi}_{N})}}{1+p\left(-1+\e^{-K_N(t,\boldsymbol{\xi}_{N})}\right)} \right) &\leq &\displaystyle \left| f_{P_0}\left( \frac{p \e^{-K_N(t,\boldsymbol{\xi}_{N})}}{1+p\left(-1+\e^{-K_N(t,\boldsymbol{\xi}_{N})}\right)} \right) -f_{P_0}(p_{0,1})\right|+ f_{P_0}(p_{0,1}) \\
\\
\displaystyle &\leq&\displaystyle  L_{f_{P_0}} \left|  \frac{p \e^{-K_N(t,\boldsymbol{\xi}_{N})}}{1+p\left(-1+\e^{-K_N(t,\boldsymbol{\xi}_{N})}\right)}-p_{0,1}\right|+F_0\\
\\
\displaystyle & \leq & \displaystyle   L_{f_{P_0}} \left(\left|  \frac{p \e^{-K_N(t,\boldsymbol{\xi}_{N})}}{1+p\left(-1+\e^{-K_N(t,\boldsymbol{\xi}_{N})}\right)}\right|+p_{0,1}\right)+F_0\\
\\
& <&\displaystyle   L_{f_{P_0}}(1+p_{0,1})+F_0.
\end{array}
\]
Observe that in the last inequality we have used Remark \ref{P_N_entre_0y1}.

Now, we will bound term (2). With this end the mean value theorem (MVT) will be applied. Let $p\in  \mathcal{J}$ be arbitrary but fixed, being $ \mathcal{J}\subset ]0,1[$ bounded and let us define the auxiliary function 
\[
g(z)=\frac{\e^{-z}}{(1+p(-1+\e^{-z}))^2}, \quad z\in \mathbb{R}.
\]

\begin{rmk}\label{rmkpolo}
We prove that $1+p(-1+\e^{-z})\neq 0$ for all $z\in \mathbb{R}$. Let us reasoning by contradiction. Assume that $1+p(-1+\e^{-z})= 0$.  The case $z=0$ obviously leads to contradiction. $1+p(-1+\e^{-z})= 0$ if and only if $p=\frac{\e^{z}}{\e^{z}-1}$. But, if $z> 0$ then $\frac{\e^{z}}{\e^{z}-1}\geq 1$, while if $z<0$ then $\frac{\e^{z}}{\e^{z}-1}< 0$. So, both cases  lead to contradiction too, since $p\in ]0,1[$.
\end{rmk}

Notice that, by Remark \ref{rmkpolo} the function $g(z)$ is well-defined. Its derivative $g'(z)$ is given by
\[
g'(z)=-\frac{\e^{z}(\e^{z}(-1+p)+p)}{(\e^{z}(-1+p)-p)^3}, \quad z\in \mathbb{R}.
\]
The function $g'(z)$ is bounded for every $ p\in \mathcal{J}$, being $\mathcal{J}\subset ]0,1[$ bounded:
\begin{itemize}
\item By Remark \ref{rmkpolo} we can assure that $(\e^{z}(-1+p)-p)\neq 0$ (just  multiplying $1+p(-1+\e^{-z})\neq 0$ by $-\e^{z}$). Thus $g'(z)$ is well-defined for all $z\in \mathbb{R}$.
\item Moreover, using L'H\^{o}pital rule, it is easy to check that: $\displaystyle \lim_{z\to \pm \infty} g'(z)=0$. 
\end{itemize}
Therefore, $\forall z\in \mathbb{R}$, $\exists C_g(p)>0$ such that $\vert g'(z) \vert \leq C_g(p)\leq \sup_{p\in \mathcal{J}}{\{C_g(p)\}}=C_g$. By applying the MVT there exists $\delta_g \in]\min\{K_N(t,\boldsymbol{\xi}_{N}),K_M(t,\boldsymbol{\xi}_{M})\},\max\{K_N(t,\boldsymbol{\xi}_{N}),K_M(t,\boldsymbol{\xi}_{M})\}[$ such that
\[
\begin{array}{lcl}
\displaystyle \left| \frac{\e^{-K_N(t,\boldsymbol{\xi}_{N})}}{\left(1+p\left(-1+\e^{-K_N(t,\boldsymbol{\xi}_{N})}\right)\right)^2}-\frac{\e^{-K_M(t,\boldsymbol{\xi}_{M})}}{\left(1+p\left(-1+\e^{-K_M(t,\boldsymbol{\xi}_{M})}\right)\right)^2}\right|&=&\vert g'(\delta_g) \vert \vert K_N(t,\boldsymbol{\xi}_{N})-K_M(t,\boldsymbol{\xi}_{M})\vert\\
\\
& \leq& C_g\vert K_N(t,\boldsymbol{\xi}_{N})-K_M(t,\boldsymbol{\xi}_{M})\vert. 
\end{array}
\]
To bound the expression (3), the same argument exhibited to obtain the bound (2) can be applied. In this case, we apply the MVT to the auxiliary function
\[
h(z)=\frac{p\e^{-z}}{1+p(-1+\e^{-z})},\quad z\in \mathbb{R},
\]
whose derivative is also bounded, $|h'(z)|\leq C_h$. Therefore, by Hypothesis \textbf{H3} and the MVT, it is guaranteed the existence of $ \delta_h \in]\min\{K_N(t,\boldsymbol{\xi}_{N}),K_M(t,\boldsymbol{\xi}_{M})\},\max\{K_N(t,\boldsymbol{\xi}_{N}),K_M(t,\boldsymbol{\xi}_{M})\}[,$ such that
\[
\begin{array}{l}
\displaystyle \left|  f_{P_0}\left( \frac{p \e^{-K_N(t,\boldsymbol{\xi}_{N})}}{1+p\left(-1+\e^{-K_N(t,\boldsymbol{\xi}_{N})}\right)}\right) -f_{P_0}\left( \frac{p \e^{-K_M(t,\boldsymbol{\xi}_{M})}}{1+p\left(-1+\e^{-K_M(t,\boldsymbol{\xi}_{M})}\right)}  \right)\right|
\\
\\
\displaystyle  \leq L_{f_{P_0}} \left| \frac{p \e^{-K_N(t,\boldsymbol{\xi}_{N})}}{1+p\left(-1+\e^{-K_N(t,\boldsymbol{\xi}_{N})}\right)} - \frac{p \e^{-K_M(t,\boldsymbol{\xi}_{M})}}{1+p\left(-1+\e^{-K_M(t,\boldsymbol{\xi}_{M})}\right)} \right|\\
\\
=L_{f_{P_0}}\vert h'(\delta_h) \vert  \vert K_N(t,\boldsymbol{\xi}_{N})-K_M(t,\boldsymbol{\xi}_{M})\vert \leq L_{f_{P_0}} C_h\vert K_N(t,\boldsymbol{\xi}_{N})-K_M(t,\boldsymbol{\xi}_{M})\vert. 
\end{array}
\]
Finally to bound the term (4), as 
\[
0<\frac{p\e^{-K_M(t,\boldsymbol{\xi}_{M})}}{1+p(-1+\e^{-K_M(t,\boldsymbol{\xi}_{M})})}<1  \quad \text{and} \quad \hat{p}< p,
\]
(recall that $\hat{p}<p$ for all $p\in \mathcal{J}$ bounded), then
\[
\frac{\e^{-K_M(t,\boldsymbol{\xi}_{M})}}{\left(1+p\left(-1+\e^{-K_M(t,\boldsymbol{\xi}_{M})}\right)\right)^2}=\left(\frac{p\e^{-K_M(t,\boldsymbol{\xi}_{M})}}{1+p(-1+\e^{-K_M(t,\boldsymbol{\xi}_{M})})}\right)^2\frac{\e^{K_M(t,\boldsymbol{\xi}_{M})}}{p^2}<\frac{\e^{K_M(t,\boldsymbol{\xi}_{M})}}{\hat{p}^2}.
\]

\textbf{Step (III):} In this step of the proof, we assume the following hypothesis

\[
\textbf{H4}:
\begin{array}{c}
A(t,\omega) \,\, \text{admits a Karhunen-Lo\`{e}ve expansion of type \eqref{KLEX} such that:} \\
\exists\, C>0: \mathbb{E}
\left[
\e^{2K_N(t,\boldsymbol{\xi}_N(\omega))}
\right]\leq C,\,\,
\text{for all positive integer}\,\, N.
\end{array}
\]
Then, denoting
\[
\alpha=\left(\left(L_{f_{P_0}}(1+p_{0,1})+F_0 \right)  C_g+\frac{L_{f_{P_0}}C_h}{\hat{p}^2}C^{1/2}\right)>0,
\]
the right part of the inequality in the Step (III) is obtained. 

Summarizing, under hypotheses \textbf{H1}--\textbf{H4} we have proved that
\[
\left| f_1^N(p,t)-f_1^M(p,t) \right| < \alpha\,\mathbb{E}\left[\vert K_N(t,\boldsymbol{\xi}_{N}(\omega))-K_M(t,\boldsymbol{\xi}_{M}(\omega))\vert^2 \right]^{1/2}.
\]
Finally, following the same argument than in \cite{KL1}, that is using Cauchy--Schwarz inequality for integrals, one gets
\[
 \mathbb{E}
  \left[
\left|
K_N(t,\boldsymbol{\xi}_{N}(\omega))- K_M(t,\boldsymbol{\xi}_{M}(\omega))
\right| ^2
\right] 
\leq
(T-t_0)
\left(
\left\|
A_N(t,\omega)-A_M(t,\omega)
\right\|_{\mathrm{L}^2(\Omega, \mathrm{L}^2([t_0,T]))}
\right)^2,
\]
and as a consequence of mean square convergence of KLE to $A(t,\omega)$, one obtains
\[
\left| f_1^N(p,t)-f_1^M(p,t) \right| < \alpha\,\sqrt{T-t_0}
\left(
\left\|
A_N(t,\omega)-A_M(t,\omega)
\right\|_{\mathrm{L}^2(\Omega, \mathrm{L}^2([t_0,T]))}
\right)\xrightarrow[N,M\to +\infty]{}  0.
\]
Summarizing, the following result has been established:
\begin{thm}\label{main_result}
Under hypotheses \textbf{H1}--\textbf{H4}, the sequence $\{f^N_1(p, t) : N \geq 1\}$ of 1-PDFs, defined by \eqref{1PDF_K_N}, converges for every $(p, t)\in \mathcal{J} \times [t_0,T]$, for all $\mathcal{J} \subset \mathbb{R}$ bounded, to the exact 1-PDF, $f_1(p, t)$, of the solution SP of random IVP \eqref{logistic_problem_random}.
\end{thm}

\section{Numerical examples}\label{sec_ejemplos}

This section is devoted to illustrate the theoretical findings previously obtained through three  numerical examples. In these examples, we will compute approximations to the 1-PDF, $f_1(p,t)$, of the solution SP of the random IVP \eqref{logistic_problem_random} via $f_1^N(p,t)$, given in \eqref{1PDF_K_N}, for different probability distributions of the initial condition, $P_0(\omega)$, and different SPs for the diffusion coefficient, $A(t,\omega)$. In the first example, the standard Wiener process, also termed Browninan motion, will play the role of $A(t,\omega)$, since, in such a case, an exact  solution to  IVP  \eqref{logistic_problem_random} is available, and then we  can check graphic and numerically the accuracy of the  approximations, $f_1^N(p,t)$, for different orders of truncation $N$. Thus, Example~\ref{Ejemplo_Numerico_1} is a test example.  In the second and third examples, exact solutions are not available. In these two latter cases, we   illustrate  convergence of approximations of the 1-PDF by means of appropriate   graphical representations  and also calculating some  measure errors that involve two consecutive approximations, namely, $f_1^N(p,t)$ and $f_1^{N-1}(p,t)$. Afterwards,  approximations of the mean and  variance of the solution SP, $P(t,\omega)$, will be   computed in  the three examples using the following expressions
\begin{eqnarray}\label{media_varianza_aprox}
\mathbb{E}\left[\left(P_N(t,\omega)\right)\right]=\int_{-\infty}^{\infty} p f_1^N(p,t)\,\mathrm{d}p\,,
\quad
\mathbb{V}\left[\left(P_N(t,\omega)\right)\right]=\int_{-\infty}^{\infty} p^2 f_1^N(p,t)\,\mathrm{d}p-(\mathbb{E}\left[\left(P_N(t,\omega)\right)\right])^2 ,
\end{eqnarray} 
where $P_N(t,\omega)$ is defined in \eqref{P_N}. Finally,  we will assess the accuracy of  these approximations (the mean and the variance) via  appropriate error measures that will be introduced in the examples.

\begin{example}\label{Ejemplo_Numerico_1}
Let us consider the random IVP \eqref{logistic_problem_random} on the time interval $\mathcal{T}=[t_0,T]=[0,1.5]$. We choose as  initial condition $P_0(\omega)$ a  Beta RV with parameters $7$ and $10$ truncated on the interval $[0.1,0.9]$, $P_0(\omega) \sim \text{Be}_{[0.1,0.9]} (7;10)$ and, as the diffusion coefficient $A(t,\omega)$, the standard Wiener process, $W(t,\omega)$, whose mean and covariance functions are given by $\mu_W(t)=0$, $t\in \mathcal{T}$ and $c_W(s,t)=\min(s,t)$, $\forall (s,t)\in \mathcal{T}\times \mathcal{T}$, respectively. The KLE of the standard Wiener process is given by \eqref{KLEX} with  $\xi_j(\omega)$  pairwise uncorrelated standard Gaussian RVs, $\xi_j(\omega)\sim \text{N}(0;1)$, and being 
\[
\nu_j=\frac{4T^2}{(2j-1)^2 \pi^2},\quad \phi_j(t)=\sqrt{\frac{2}{T}}\sin \left(\frac{(2j-1)\pi t}{2T}\right), \quad j=1,2,\dots
\]
the corresponding eigenvalues and eigenfunctions, respectively, \cite[p.~206]{Powell-LIBRO}. We will choose $P_0(\omega)$ so that is independent of the random vector $\boldsymbol{\xi}_{N}=(\xi_1(\omega),\ldots,\xi_N(\omega))$, for $N$ arbitrary but fixed.

Let us check that hypotheses \textbf{H1}--\textbf{H4} of Th.~\ref{main_result} hold. Since $P_0(\omega) \sim \text{Be}_{[0.1,0.9]} (7;10)$, then first part of hypothesis \textbf{H1} is clearly satisfied taking $p_{0,1}=0.1$ and $p_{0,2}=0.9$. To check the second part, it is enough to observe that $\int_{t_0}^{T} \mathbb{E}[(W(t,\omega))^2]\, \mathrm{d} t =\int_{t_0}^{T}  t\, \mathrm{d} t=\frac{T^2}{2}-\frac{t_0^2}{2}=1.125<\infty $.  The hypothesis \textbf{H2} holds because the choice we have made for initial condition, $P_0(\omega)$, and random vector $\boldsymbol{\xi}_{N}(\omega)$. It is straightforward to check that the first derivative of the PDF of  $P_0(\omega) \sim \text{Be}_{[0.1,0.9]} (7;10)$ is bounded over the domain $[0.1,0.9]$, thus $f_{P_0}(p_0)$ is Lipschitz on $[0.1,0.9]$. This justifies hypothesis \textbf{H3}. Finally, following the same reasoning exhibited in \cite[Remark 2]{KL1}, it is checked that hypothesis \textbf{H4} fulfils.

Thus according to \eqref{1PDF_K_N}, the  1-PDF of the approximate solution SP, $P_N(t,\omega)$,   is given by
\begin{equation}\label{1-PDF-Brownian}
f_1^N(p,t)=\int_{\mathbb{R}^N} f_{P_0}\left(\frac{p \e^{-\sum_{j=1}^Nh_j(t)\xi_j}}{1+p(-1+\e^{-\sum_{j=1}^Nh_j(t)\xi_j})}\right)\frac{\prod_{j=1}^N f_{\xi_j}(\xi_j) \e^{-h_j(t)\xi_j}}{\left(1+p\left(-1+\e^{-\sum_{j=1}^Nh_j(t)\xi_j}\right)\right)^2} \mathrm{d} \xi_N \cdots \mathrm{d}\xi_1,
\end{equation}
where $f_{\xi_j}(\xi_j)$ denotes the PDF of $\xi_j(\omega)$ for each $j:1 \leq j \leq N$ and
\begin{equation}\label{hj}
h_j(t)=\left(\frac{2T}{(2j-1)\pi}\right)^2\sqrt{\frac{2}{T}}\left(1-\cos\left(\frac{(2j-1)\pi t}{2T}\right)\right).
\end{equation}

So far we have obtained the approximations $f_1^N(p,t)$ to the exact 1-PDF, $f_1(p,t)$, of the solution SP of the random IVP \eqref{logistic_problem_random}. Now, we will determine an explicit expression to $f_1(p,t)$. To this goal, let us recall that $\int_{0}^t W(s,\omega)\mathrm{d}s \sim \text{N}\left(0;\sqrt{\frac{t^3}{3}}\right)$, \cite[p.~105]{Soong}, hence $\int_{0}^t W(s,\omega)\mathrm{d}s \stackrel[]{\text{d}}{=} \sqrt{\frac{t^3}{3}}Z(\omega)$, $Z(\omega)\sim \text{N}(0;1)$. Taking into account that the solution SP of the random IVP \eqref{logistic_problem_random} is expressed in terms of this stochastic integral (see \eqref{logistic_solution} with $A(t,\omega)=W(t)$) and, by applying the  RVT technique (see Th.~\ref{RVT}), it is straightforwardly to check that
\begin{equation}\label{1-PDF-Brownian-Exact}
f_1(p,t)=\int_{\mathbb{R}} f_{P_0} \left(\frac{1}{1+\exp\left(\sqrt{\frac{t^3}{3}}z\left(-1+\frac{1}{p}\right)\right)}\right)f_Z(z) \frac{\exp\left(\sqrt{\frac{t^3}{3}}z\right)}{\left(p+\exp\left(\sqrt{\frac{t^3}{3}}z\right)(1-p)\right)^2} \mathrm{d}z
\end{equation}

In Figure \ref{grafica1PDFs}, we show the 1-PDF, $f_1(p,t)$, given by \eqref{1-PDF-Brownian-Exact} of the exact solution SP \eqref{logistic_solution_random} together with the 1-PDFs, $f_1^N(p,t)$, given by \eqref{1-PDF-Brownian}--\eqref{hj} corresponding to  the approximate  solution SP \eqref{P_N}  with $N\in\{1,2\}$. We can observe that these approximations  $f_1^N(p,t)$ clearly converge to $f_1(p,t)$, even for small values of the order of truncation $N$.   For the sake of clarity, in Figure \ref{grafica1PDFs-t} we have plotted the exact 1-PDF, $f_1(p,t)$, and the approximate 1-PDFs, $f_1^N(p,t)$, for $N\in\{1,2,3\}$  at different time instants,  $t\in\{0.50,0.75,1.00,1.50\}$. Again, we can observe fast convergence of $f_1^N(p,t)$ to $f_1(p,t)$. In order to better assess this convergence, in Table \ref{tabla_errores} we have collected the total probabilistic error defined in \eqref{error-PDF}. From these figures we can observe that for $t$ fixed, the error $e_N^{\text{PDF}}(t)$ decreases as $N$ increases, as expected. 
\begin{equation}\label{error-PDF}
e_N^{\text{PDF}}(t)=\int_{0}^{1} \left|  f_1(p,t)-f_1^N(p,t) \right| \mathrm{d}p.
\end{equation}

\begin{figure}[htp]
\begin{center}
\includegraphics[width=0.32\textwidth]{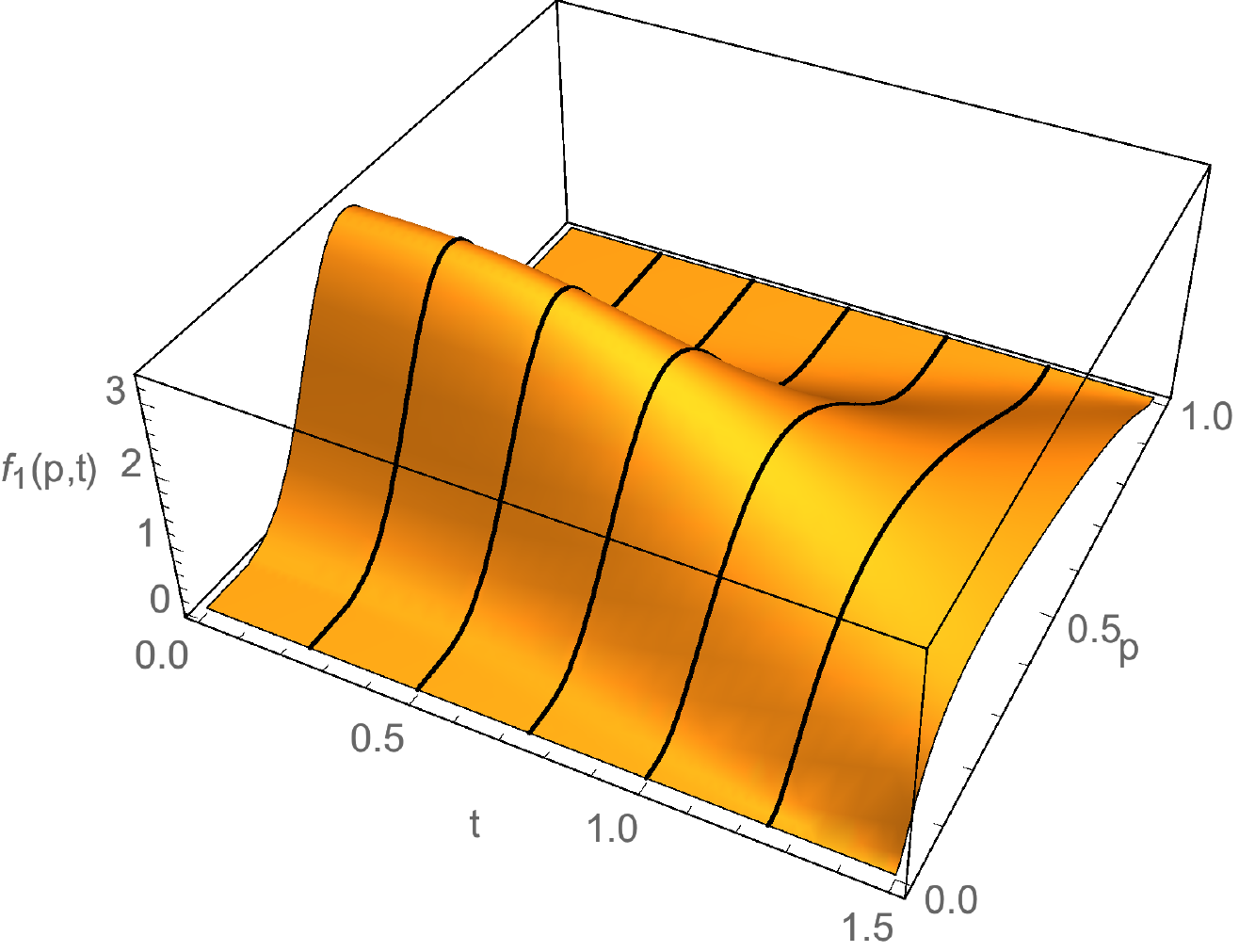}
\includegraphics[width=0.32\textwidth]{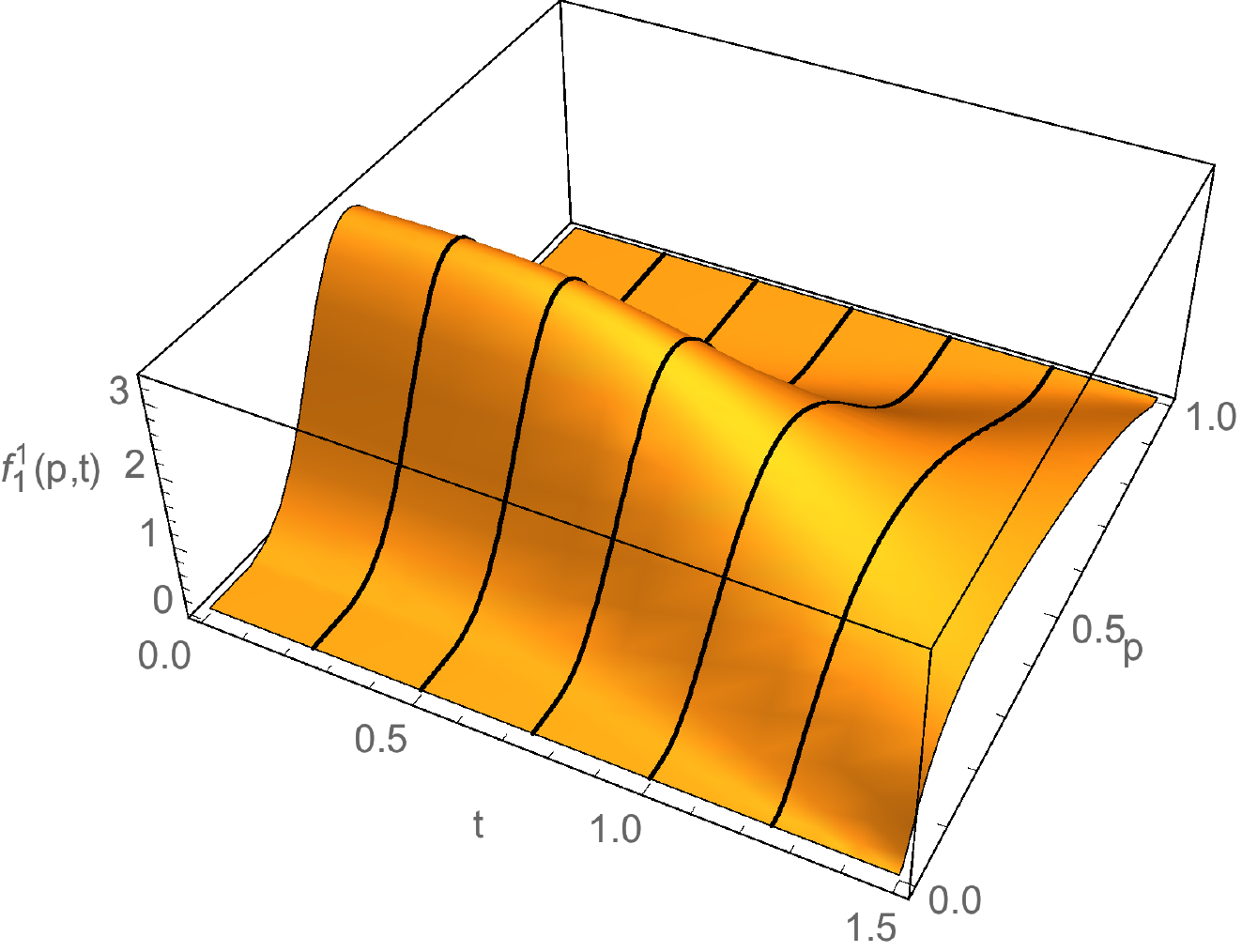}
\includegraphics[width=0.32\textwidth]{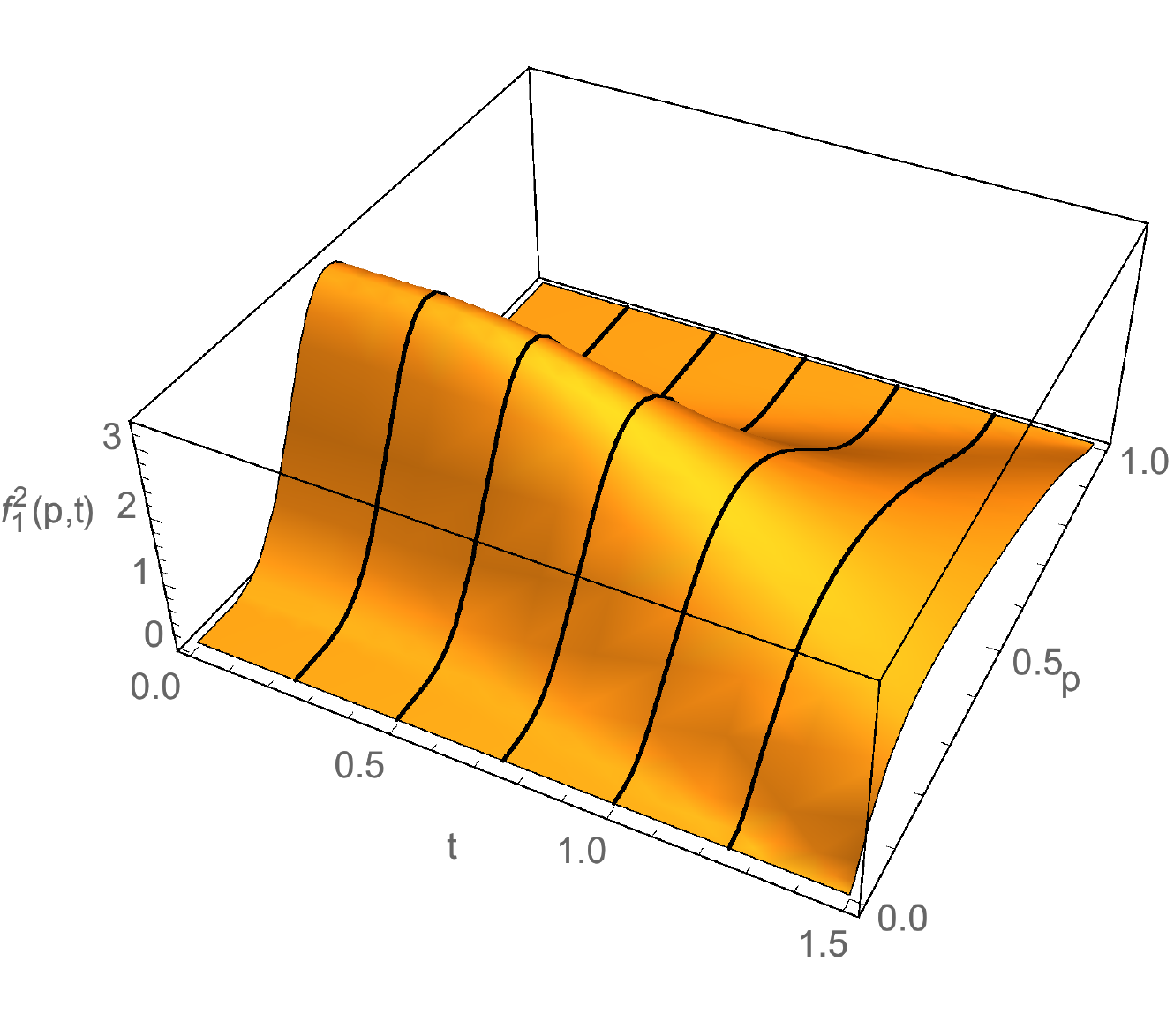}
\end{center}\caption{Example \ref{Ejemplo_Numerico_1}. Left: 1-PDF, $f_1(p,t)$, of the exact solution SP given by \eqref{1-PDF-Brownian-Exact}. Center: 1-PDF, $f_1^N(p,t)$,  of the approximate solution given by \eqref{1-PDF-Brownian}--\eqref{hj} with $N=1$. Right: 1-PDF, $f_1^N(p,t)$,  of the approximate solution given by \eqref{1-PDF-Brownian}--\eqref{hj} with $N=2$.} \label{grafica1PDFs}
\end{figure}

\begin{figure}
\begin{center}
\includegraphics[width=0.49\textwidth]{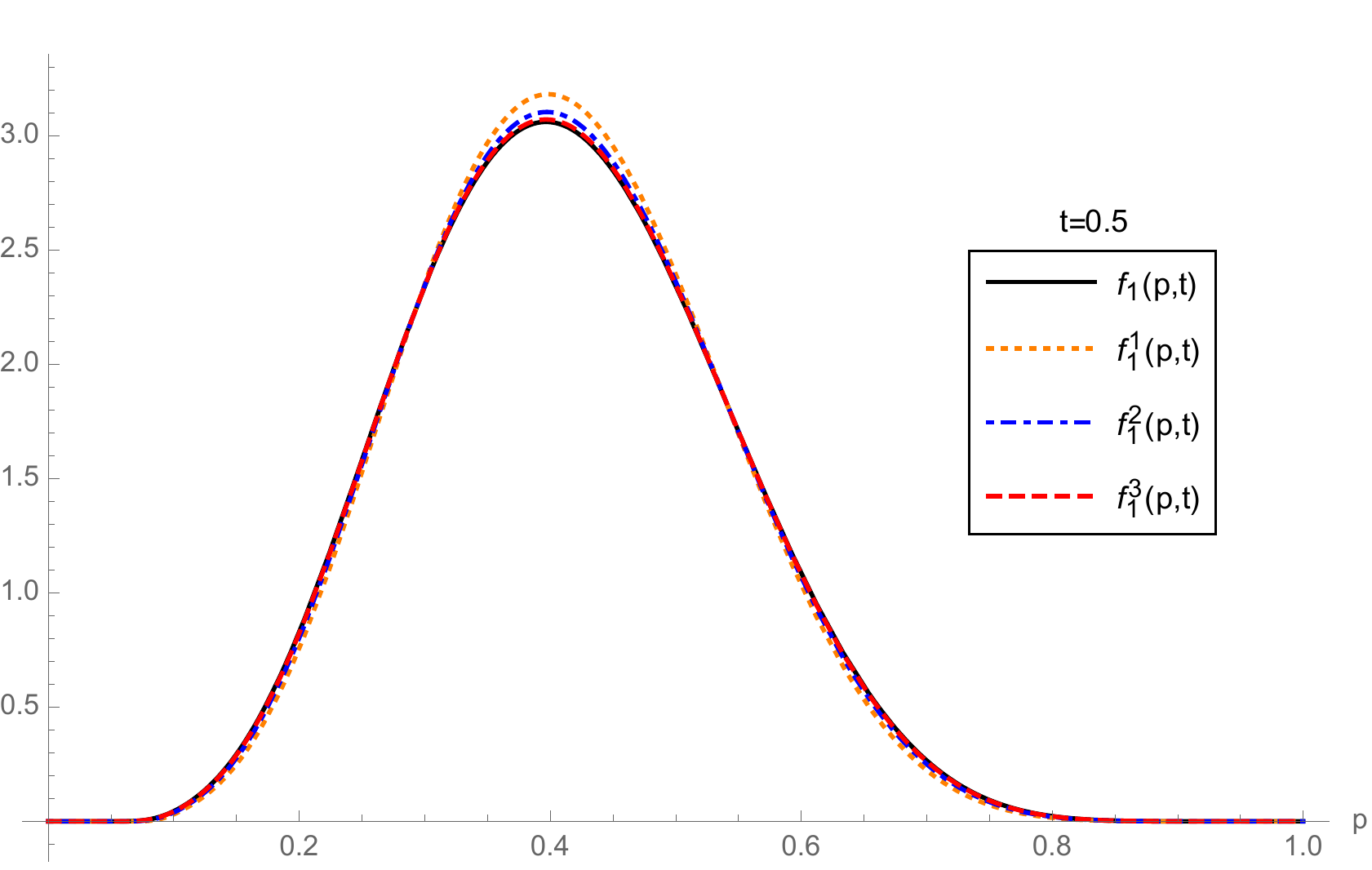}
\includegraphics[width=0.49\textwidth]{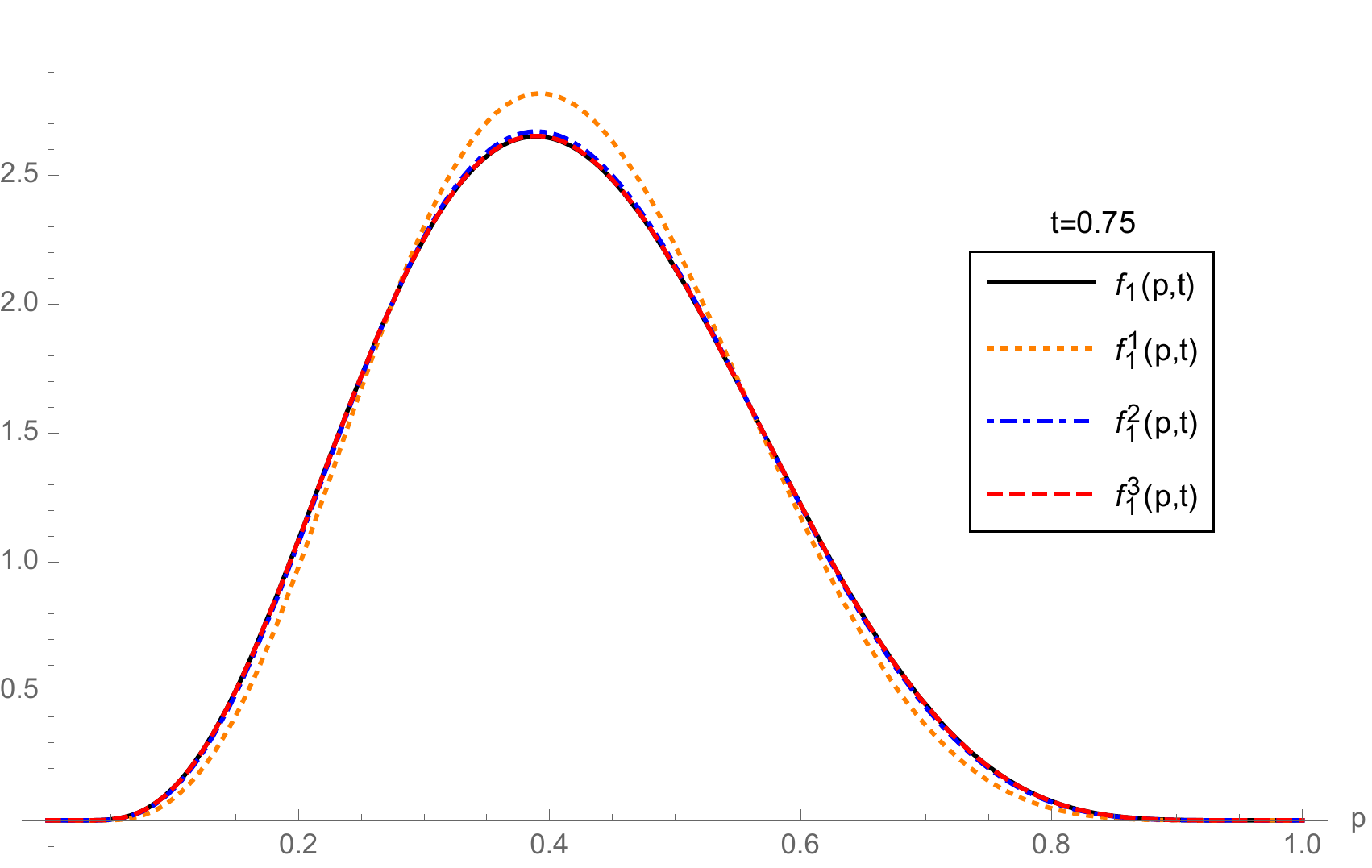}\\
\includegraphics[width=0.49\textwidth]{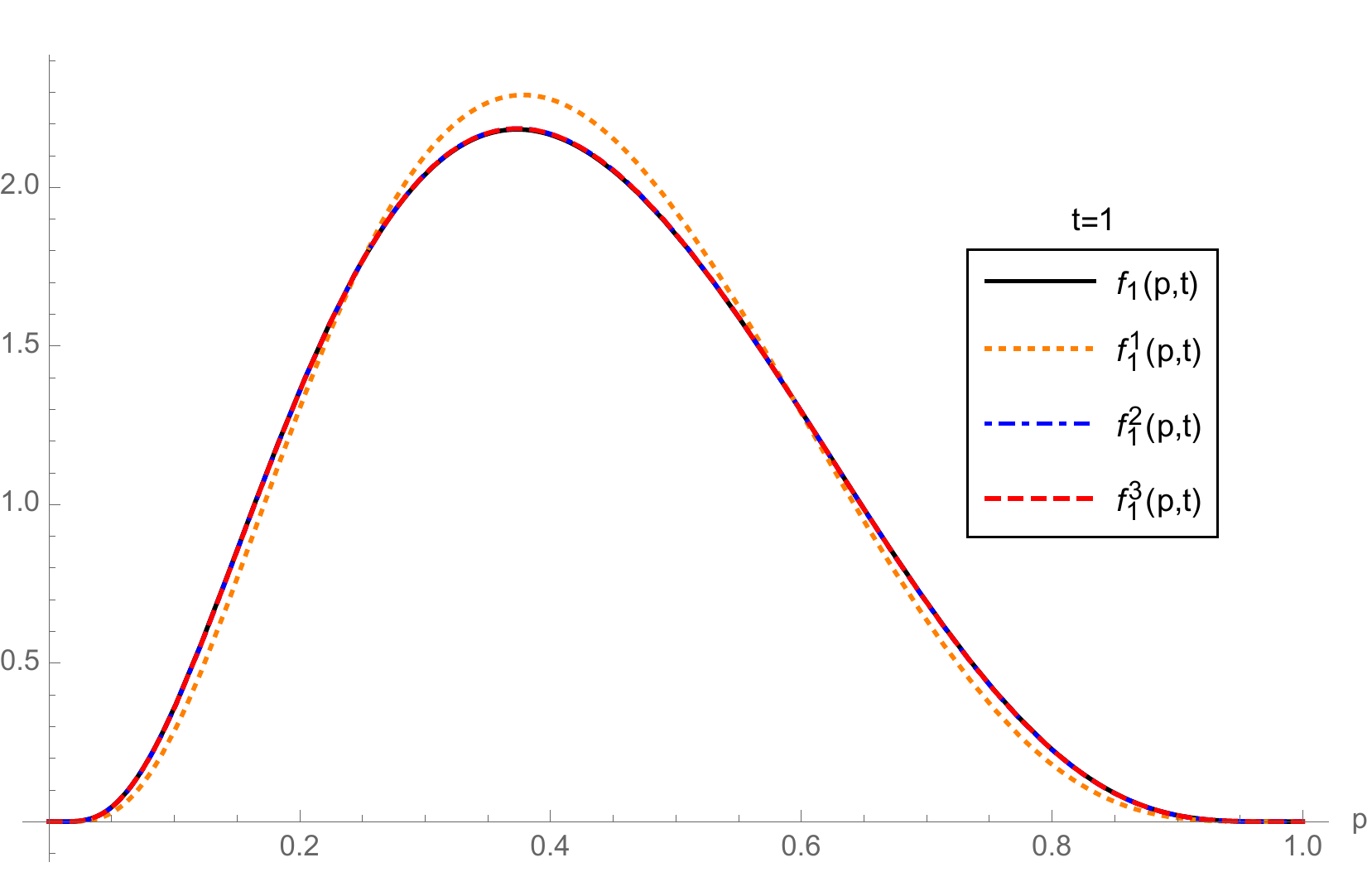}
\includegraphics[width=0.49\textwidth]{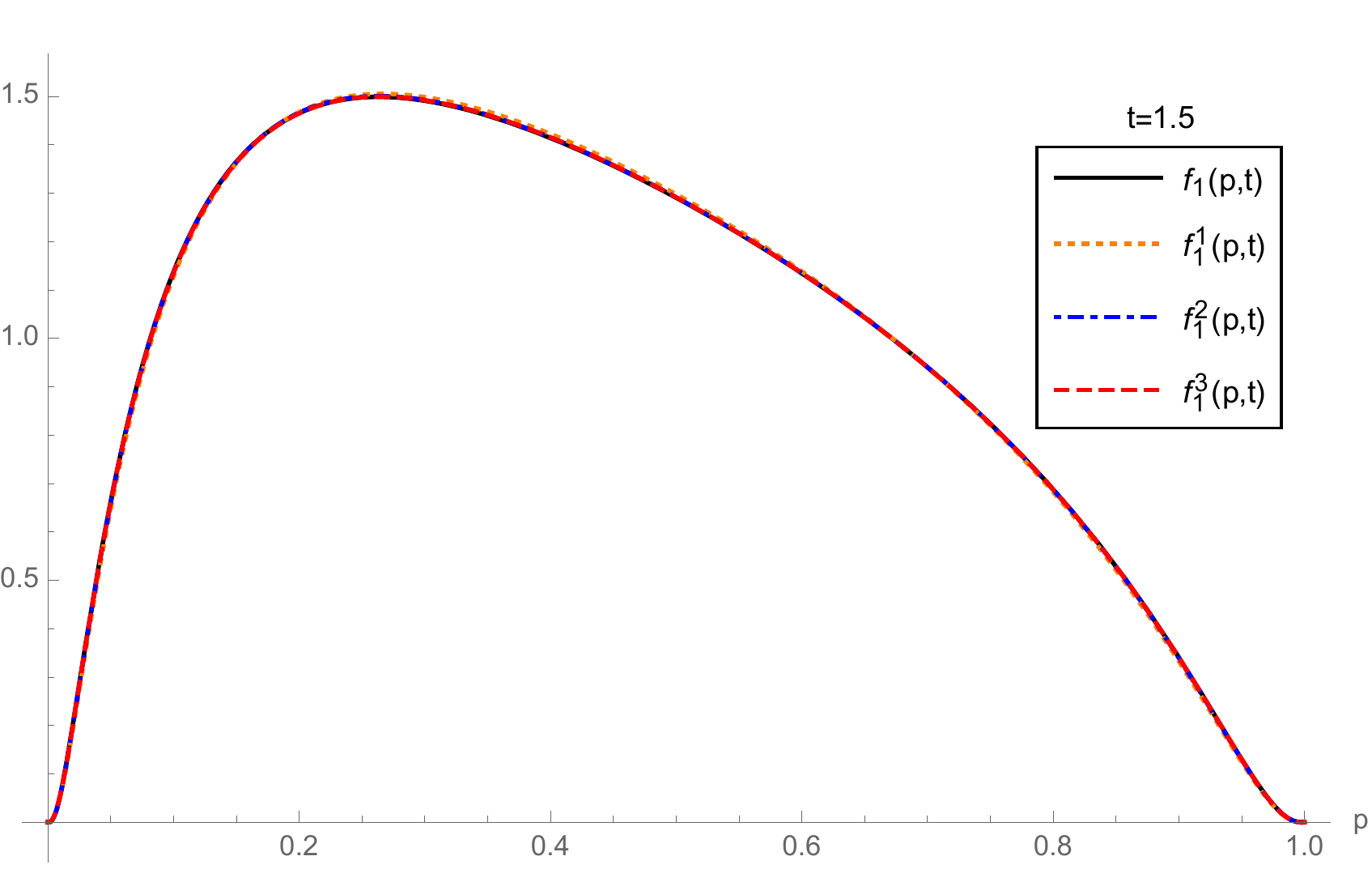}
\end{center}\caption{Example \ref{Ejemplo_Numerico_1}. Plots of the 1-PDF, $f_1(p,t)$, of the exact solution SP given by \eqref{1-PDF-Brownian-Exact} and the  truncations, $f_1^N(x,t)$, given by \eqref{1-PDF-Brownian}--\eqref{hj} with $N=1,2,3$ for different values of  $t$. Left-up:  $t=0.50$, Rigth-up: $t=0.75$. Left-down: $t=1.00$. Right-down: $t=1.50$.}\label{grafica1PDFs-t}
\end{figure}

\begin{table}[htp]
\begin{center}
\begin{tabular}{*{4}{|c}|}
\hline
$e_N^{\text{PDF}}(t)$  & $N=1$ & $N=2$ & $N=3$    \\
\hline
$t=0.50$ &   0.037418 & 0.013544 & 0.003149 \\
\hline 
$t=0.75$ &   0.059518 & 0.006964 & 0.000652 \\
\hline 
$t=1.00$ &  0.048595 & 0.001153 & 0.000987 \\
\hline 
$t=1.50$ & 0.005737 & 0.000789 & 0.000648 \\
\hline 
\end{tabular}
\end{center}
\caption{Error measure $e_N^{\text{PDF}}(t)$, defined by \eqref{error-PDF}, for different time instants, $t\in \{ 0.50,0.75,1.00,1.50 \}$, and truncation orders $N\in\{1,2,3\} $, in the context of Example \ref{Ejemplo_Numerico_1}.}
\label{tabla_errores}
\end{table}

We complete the numerical study by computing the approximations of the mean and the variance functions using \eqref{media_varianza_aprox} and \eqref{1-PDF-Brownian}--\eqref{hj} with $N\in\{1,2,3,4\}$. In Figure \ref{grafica-Mean-Variance}, we have plotted these approximations together with the exact mean and variance functions obtained via \eqref{moments}, with $k=1,2$, and \eqref{varianza_formula}, where $f_1(p,t)$ is given by \eqref{1-PDF-Brownian-Exact}. To assess the quality of these approximations, we have computed  the total  error of the mean and the variance using the following expressions
\begin{equation}\label{errorMEDIA_VARIANZA}
e_N^{\mathbb{E}}=\int_{t_0}^{T} \left| \mathbb{E}[P(t,\omega)]-\mathbb{E}[P_N(t,\omega)]\right| \mathrm{d} t,\qquad
e_N^{\mathbb{V}}=\int_{t_0}^{T} \left| \mathbb{V}[P(t,\omega)]-\mathbb{V}[P_N(t,\omega)]\right| \mathrm{d} t.
\end{equation}
In Table \ref{tabla_errores_moments_ejemplo_1}, we show the values of errors $e_N^{\mathbb{E}}$ and $e_N^{\mathbb{V}}$ for $N\in\{1,2,3,4\}$. From these figures we can observe that both errors decreases as $N$ increases, thus showing fully agreement with the graphical representation shown in Figure \ref{grafica-Mean-Variance}.

\begin{figure}
\begin{center}
\includegraphics[width=0.49\textwidth]{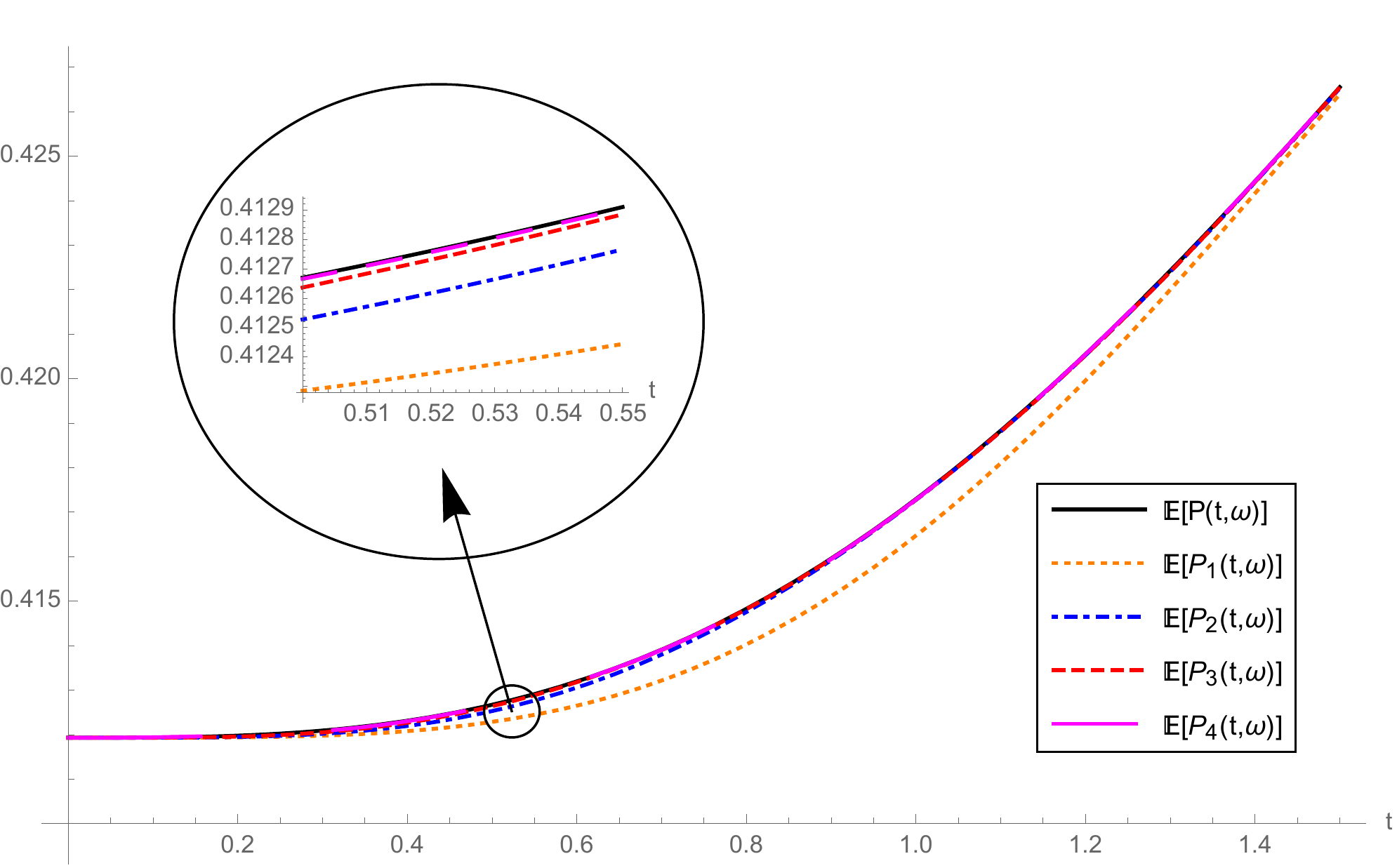}
\includegraphics[width=0.49\textwidth]{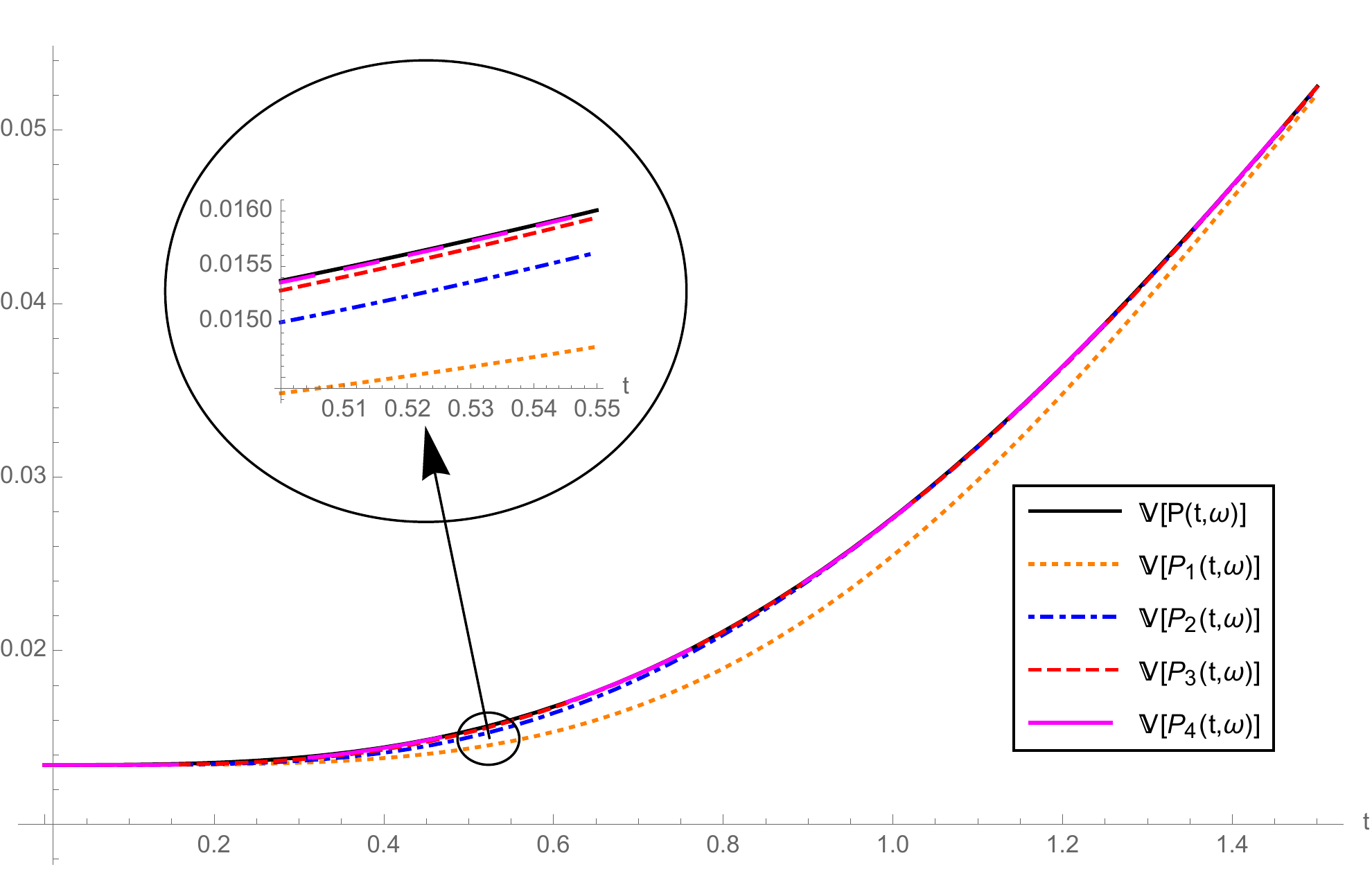}
\end{center}\caption{Example 1. Left: Exact mean, $\mathbb{E}[P(t,\omega)]$, of the solution SP and its approximations using truncations of order $N\in\{1,2,3,4\}$ ($\mathbb{E}[P_i(t,\omega)]$, $i=1,2,3,4$). Right: Exact variance ($\mathbb{V}[P(t,\omega)]$) of the solution and its approximations using truncations of order $N\in\{1,2,3,4\}$ ($\mathbb{V}[P_i(t,\omega)]$, $i=1,2,3,4$).}\label{grafica-Mean-Variance}
\end{figure}

\begin{table}[htp]
\begin{center}
\begin{tabular}{*{5}{|c}|}
\hline
$\text{Error}$  & $N=1$ & $N=2$  & $N=3$ & $N=4$  \\
\hline
Mean $e_N^{\mathbb{E}}$ &  0.000659 & 0.000085 & 0.000029 & 0.000009\\
\hline 
Variance $e_N^{\mathbb{V}}$ &   0.001756 & 0.000225 & 0.000077  & 0.000035 \\
\hline 
\end{tabular}
\end{center}
\caption{Values of errors $e_N^{\mathbb{E}}$ and $e_N^{\mathbb{V}}$ for the mean and variance, respectively, given by \eqref{errorMEDIA_VARIANZA} using different orders of  truncation $N\in\{1,2,3,4\}$, in the context of Example \ref{Ejemplo_Numerico_1}.}
\label{tabla_errores_moments_ejemplo_1}
\end{table}
\end{example}

\begin{example}\label{Ejemplo_Numerico_2}
Now we will consider the random IVP \eqref{logistic_problem_random} on the time interval $\mathcal{T}=[t_0,T]=[0,1]$. We assume that the initial condition $P_0(\omega)$ has  a truncated Exponential distribution on the interval $[0.1,0.9]$ and with  parameter $\lambda=10$, i.e. $P_0(\omega) \sim \text{Exp}_{[0.1,0.9]} (10)$. For the diffusion coefficient, $A(t,\omega)$, we will choose the Brownian Bridge \cite[p.~193--195]{Powell-LIBRO}. This SP, say $X(t,\omega)$, is defined in terms of the Wiener SP as $X(t,\omega)=W(t,\omega)-tW(1,\omega)$, having zero-mean, $\mu_X(t)=0$, and  correlation function 
\[
c_X(s,t)=\min(s,t)-st, \quad (s,t) \in \mathcal{T}\times \mathcal{T}.
\]
In \cite[p.~204]{Powell-LIBRO}, it is shown that the KLE of the Brownian Bridge is given by \eqref{KLEX} being $\xi_j(\omega)\sim \text{N}(0,1)$ pairwise uncorrelated RVs   and 
\[
 \nu_j=\frac{1}{\pi^2 j^2}, \quad \phi_j(t)=\sqrt{2} \sin (j\pi t), \quad t\in \mathcal{T}, j=1,2,\ldots.
\]
We will choose the random initial condition, $P_0(\omega)$, so that is independent of the random vector $\boldsymbol{\xi}_{N}(\omega)=(\xi_1(\omega),\ldots,\xi_N(\omega))$, for $N$ arbitrary but fixed. Analogously to Example~\ref{Ejemplo_Numerico_1}, it can be checked that hypotheses \textbf{H1}--\textbf{H4} fulfil.
Therefore, according to \eqref{1PDF_K_N}, the 1-PDF of the approximate solution SP, $P_N(t,\omega)$,  is given by
\begin{equation}\label{1PDF-ejemplo2}
f_1^N(p,t)=\int_{\mathbb{R}^N} f_{P_0}\left( \frac{p \prod_{j=1}^N \e^{-h^*_j (t) \xi_j}}{1+p\left(-1+\prod_{j=1}^N \e^{-h^*_j (t) \xi_j}\right)} \right) \frac{\prod_{j=1}^N \e^{-h^*_j (t) \xi_j} f_{\xi_j}(\xi_j)}{\left(1+p\left(-1+\prod_{j=1}^N \e^{-h^*_j (t) \xi_j}\right)\right)^2} \mathrm{d}\xi_N \dots \mathrm{d} \xi_1,
\end{equation}
where
\begin{equation}\label{hj_ejemplo2}
h^*_j(t)=\frac{\sqrt{2}}{(\pi j)^2} (1-\cos (j\pi t)).
\end{equation}

In Figure \ref{grafica1PDFs_ejemplo2}, we show the surface corresponding to the 1-PDF, $f_1^N(p,t)$, given in \eqref{1PDF-ejemplo2}--\eqref{hj_ejemplo2} for  $N=1$ and $N=2$. We can observe that both approximations are very similar, then showing convergence. For the sake of clarity, in Figure \ref{grafica1PDFs-t_ejemplo2}, we have represented the 1-PDF, $f_1^N(p,t)$, at different fixed time instants $t\in \{0.25,0.40,0.50\}$, and  increasing the order of truncation $N\in\{ 1,2,3,4\}$. From these graphical representations, we clearly observe fast convergence on the whole  domain.  To illustrate numerically   this convergence, in Table \ref{tabla_errores_ejemplo2} we show the total difference  between two consecutive  approximations, $f_1^N(p,t)$ and $f_1^{N-1}(p,t)$, at the time instants previously indicated, using the following error formula
\begin{equation}\label{error-PDF_trunca}
\hat{e}_N^{PDF}(t)=\int_{0}^{1} \left|  f_1^N(p,t)-f_1^{N-1}(p,t) \right| \mathrm{d}p, \quad N=2,3,\ldots.
\end{equation}
From figures in Table  \ref{tabla_errores_ejemplo2} we can observe that for $t$ fixed, the error $\hat{e}_N^{PDF}(t)$ decreases as $N$ increases.

\begin{figure}[htp]
\begin{center}
\includegraphics[width=0.49\textwidth]{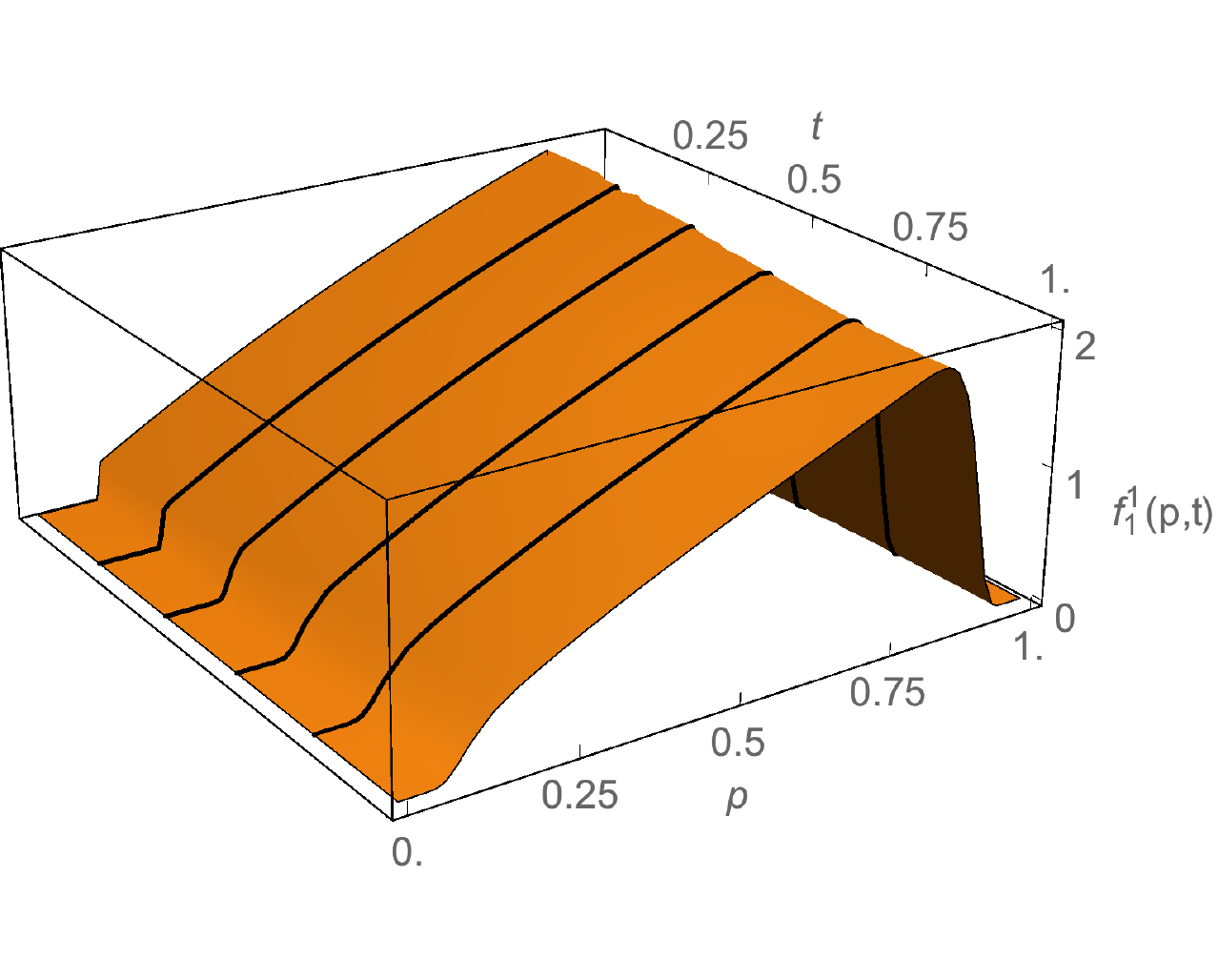}
\includegraphics[width=0.49\textwidth]{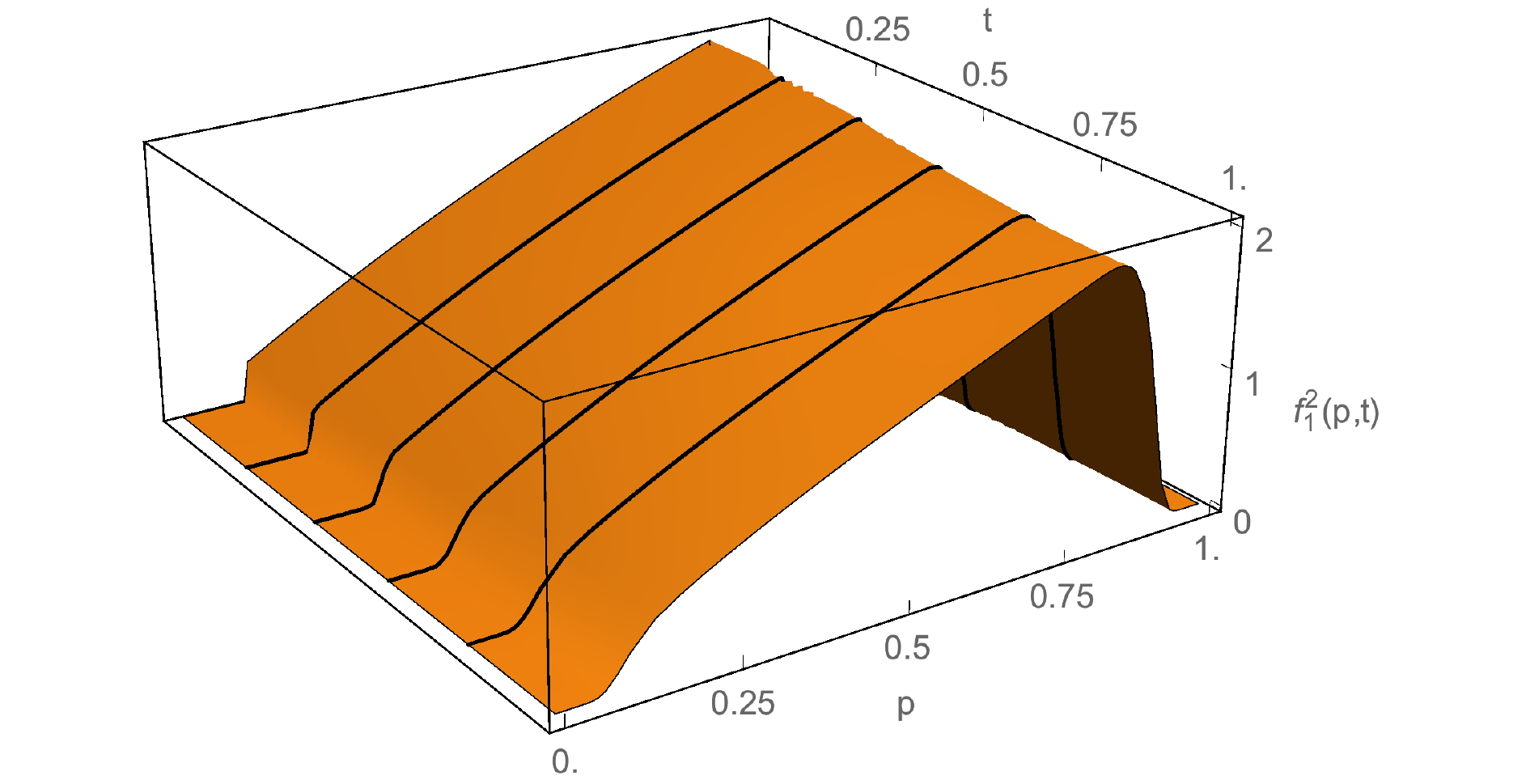}
\end{center}\caption{Example \ref{Ejemplo_Numerico_2}. Surfaces of the 1-PDF, $f_1^N(p,t)$, given in \eqref{1PDF-ejemplo2}--\eqref{hj_ejemplo2} for  $N=1$ (Left) and $N=2$ (Right).}\label{grafica1PDFs_ejemplo2}
\end{figure}

\begin{figure}
\begin{center}
\includegraphics[width=0.32\textwidth]{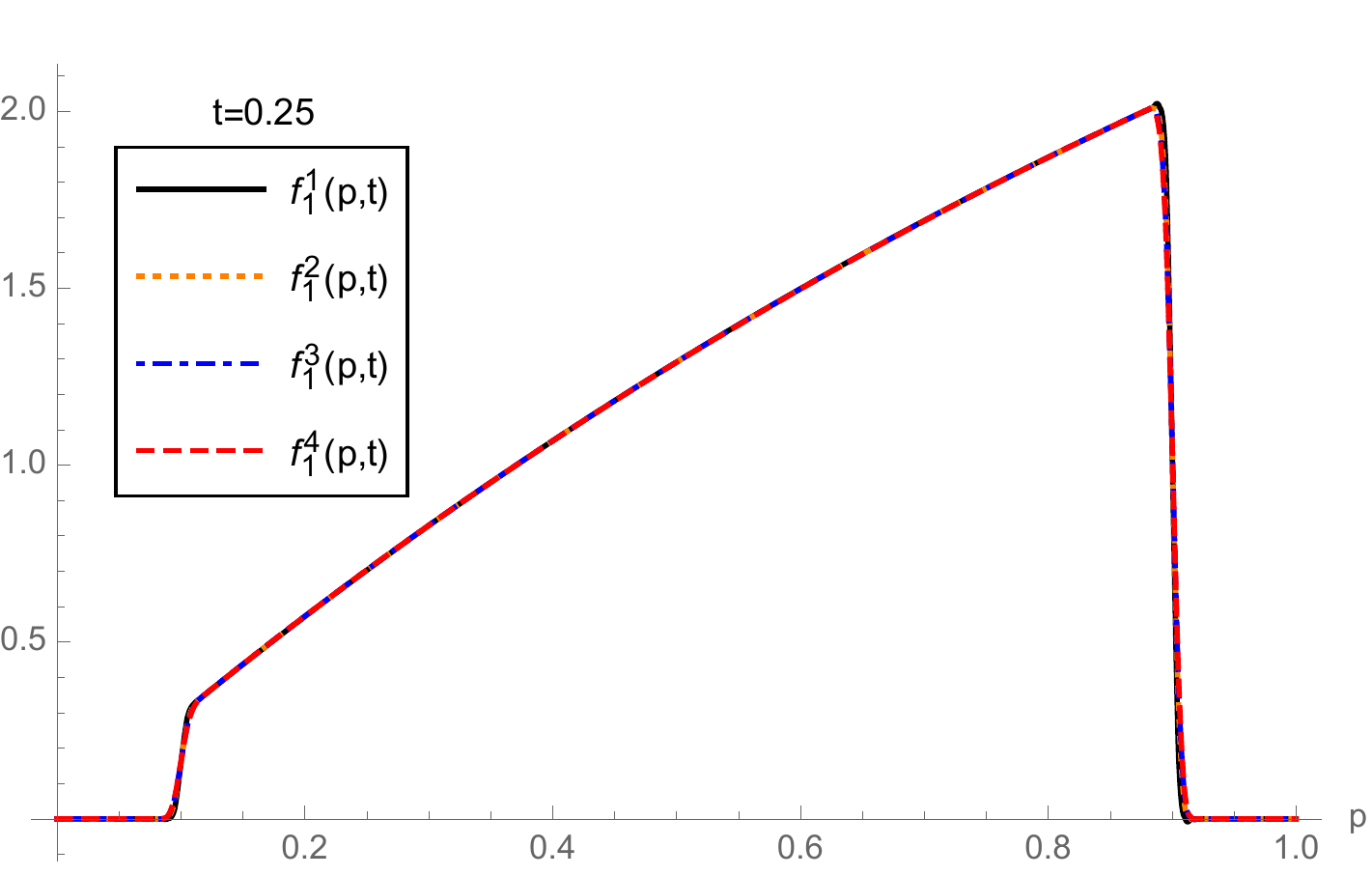}
\includegraphics[width=0.32\textwidth]{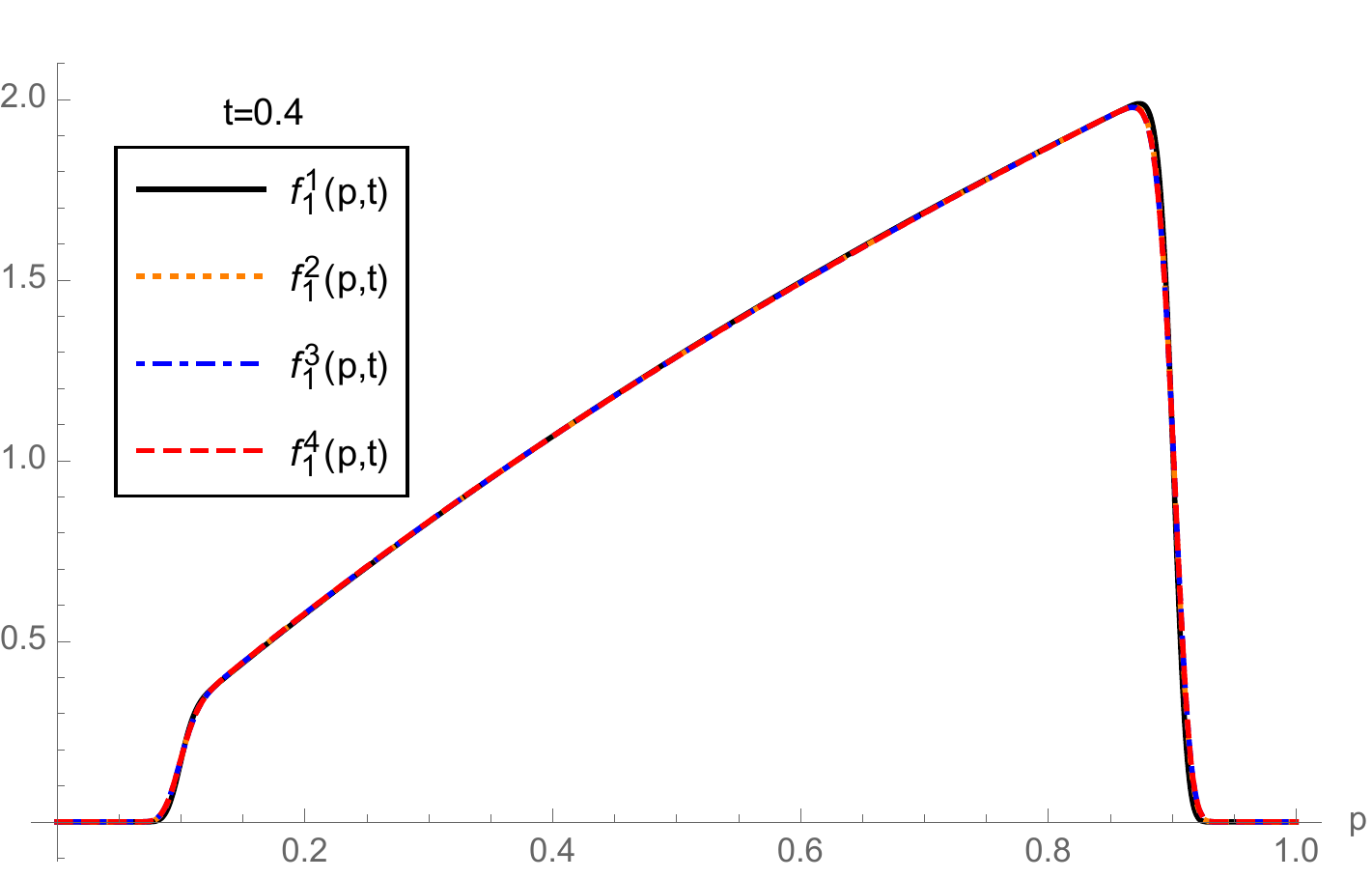}
\includegraphics[width=0.32\textwidth]{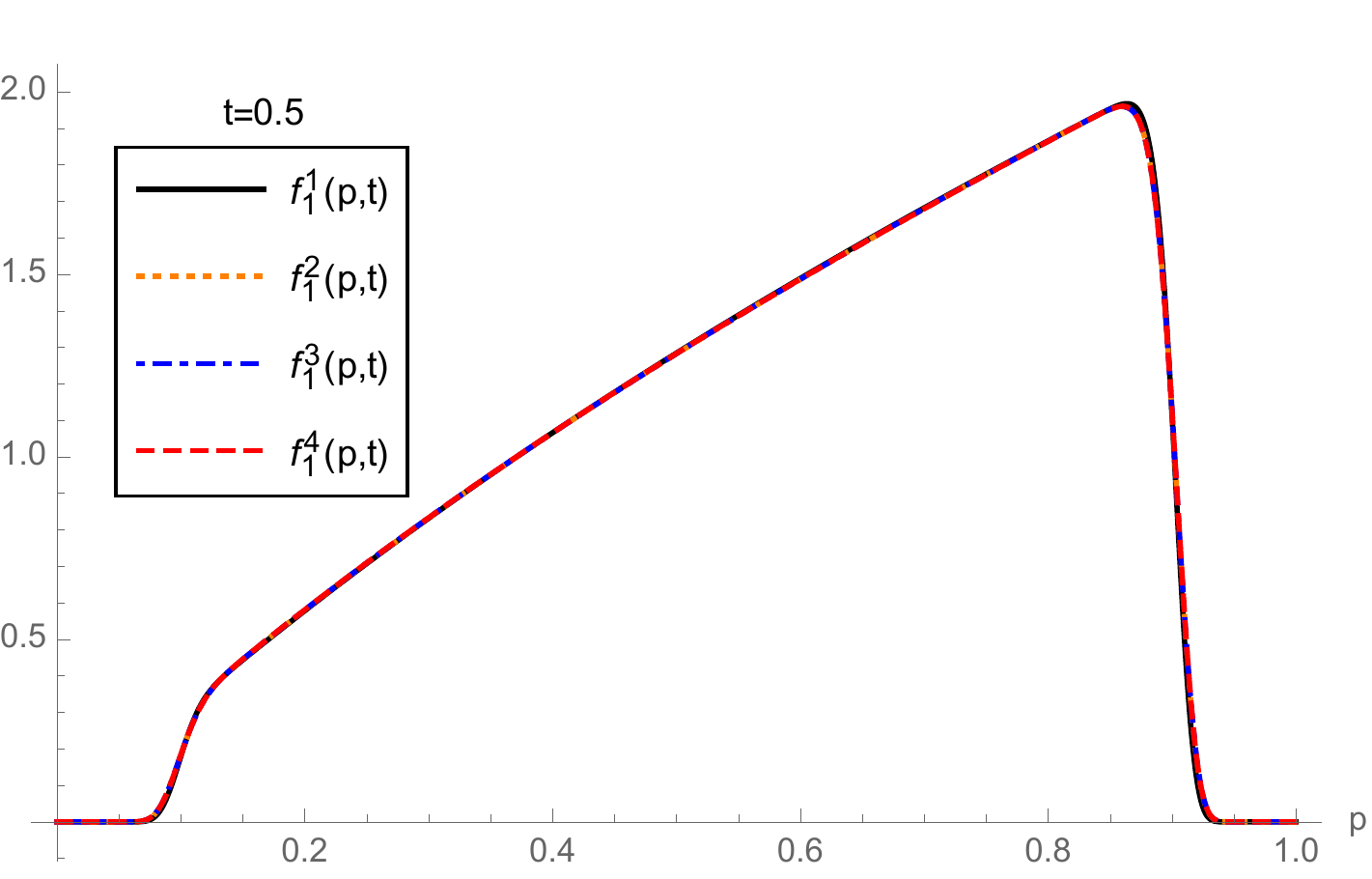}
\end{center}\caption{Example \ref{Ejemplo_Numerico_2}. Curves of the 1-PDF, $f_1^N(p,t)$, given in \eqref{1PDF-ejemplo2}--\eqref{hj_ejemplo2} at  three different time instants $t=0.25$ (Left), $t=0.40$ (Center) and $t=0.50$ (Right) using, in each case, different orders of truncations $N\in \{ 1,2,3,4\}$.}\label{grafica1PDFs-t_ejemplo2}
\end{figure}

\begin{table}[htp]
\begin{center}
\begin{tabular}{*{4}{|c}|}
\hline
$\hat{e}_N^{\text{PDF}}(t)$  & $N=2$ & $N=3$ & $N=4$    \\
\hline
$t=0.25$ &    0.002382 & 0.001275 & 0.000604\\
\hline 
$t=0.40$ &   0.004166 & 0.000746 & 0.000252 \\
\hline 
$t=0.50$ &  0.003935 & 0.000471 & 0.000306 \\
\hline 
\end{tabular}
\end{center}
\caption{Error measure $\hat{e}_N^{\text{PDF}}(t)$ defined by \eqref{error-PDF_trunca} at different time instants, $t\in \{ 0.25,0.40,0.50 \}$, and different orders of truncation  $N\in\{2,3,4\} $, in the context of Example \ref{Ejemplo_Numerico_2}.}
\label{tabla_errores_ejemplo2}
\end{table}

Finally, we take advantage of the approximations, $f_1^N(p,t)$, of the 1-PDF of the solution SP to compute approximations of the mean, $\mathbb{E}[P_N(t,\omega)]$,  and the variance, 
$\mathbb{V}[P_N(t,\omega)]$, for different orders of truncation $N\in\{1,2,3,4\}$. These approximations have been plotted in Figure \ref{grafica-Mean-Variance_ejemplo2}. In Table \ref{tabla_errores_moments}, we show the values of the following  error measure between consecutive approximations of the mean and the variance on the whole time interval.

\begin{equation}\label{errorMEDIA_VARIANZA_trunca_ej_2}
\hat{e}_N^{\mathbb{E}}=\int_{0}^{1} \left| \mathbb{E}[P_{N}(t,\omega)]-\mathbb{E}[P_{N-1}(t,\omega)]\right| \mathrm{d} t,\qquad
\hat{e}_N^{\mathbb{V}}=\int_{0}^{1} \left| \mathbb{V}[P_{N}(t,\omega)]-\mathbb{V}[P_{N-1}(t,\omega)]\right| \mathrm{d} t.
\end{equation}

\begin{figure}
\begin{center}
\includegraphics[width=0.49\textwidth]{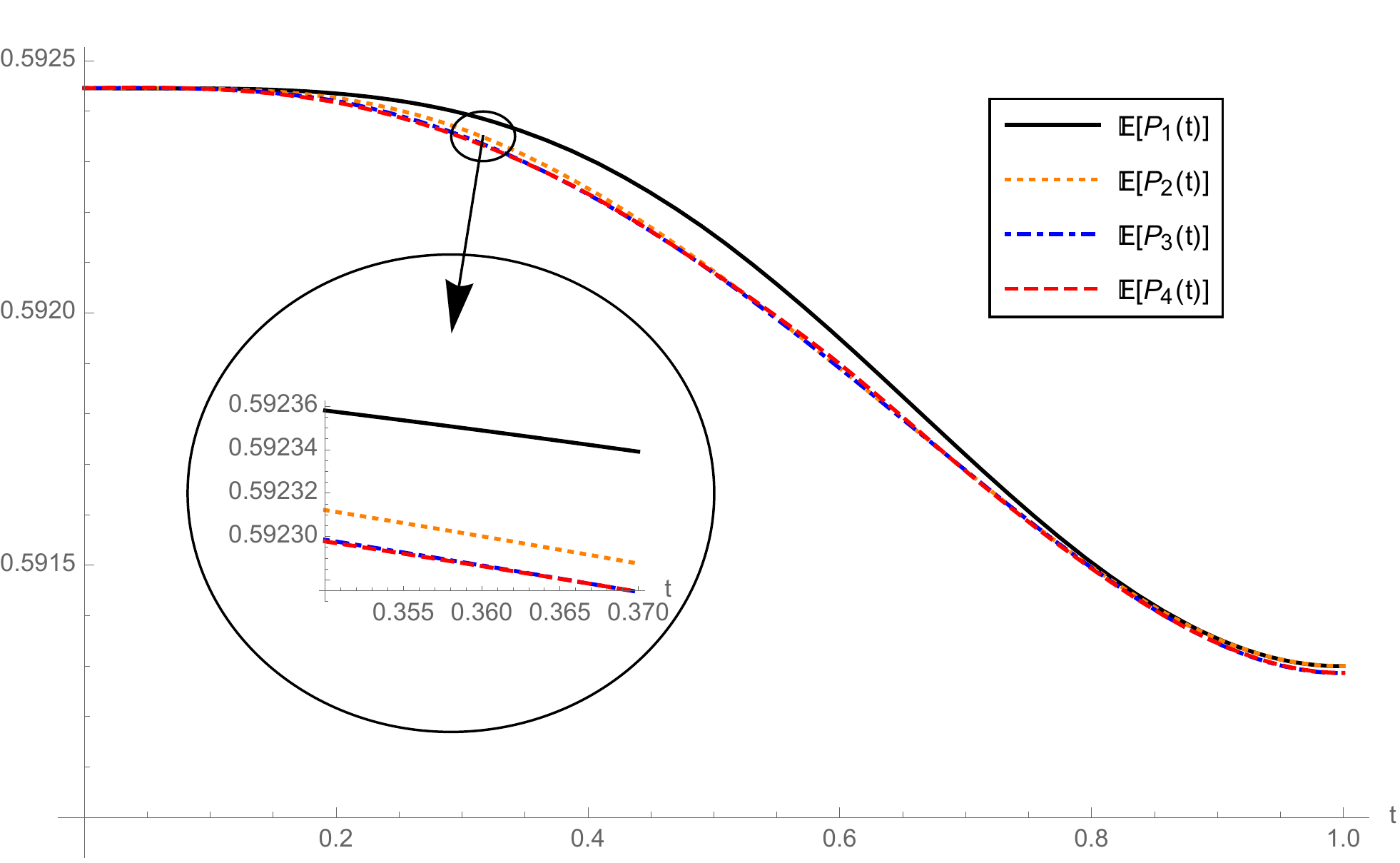}
\includegraphics[width=0.49\textwidth]{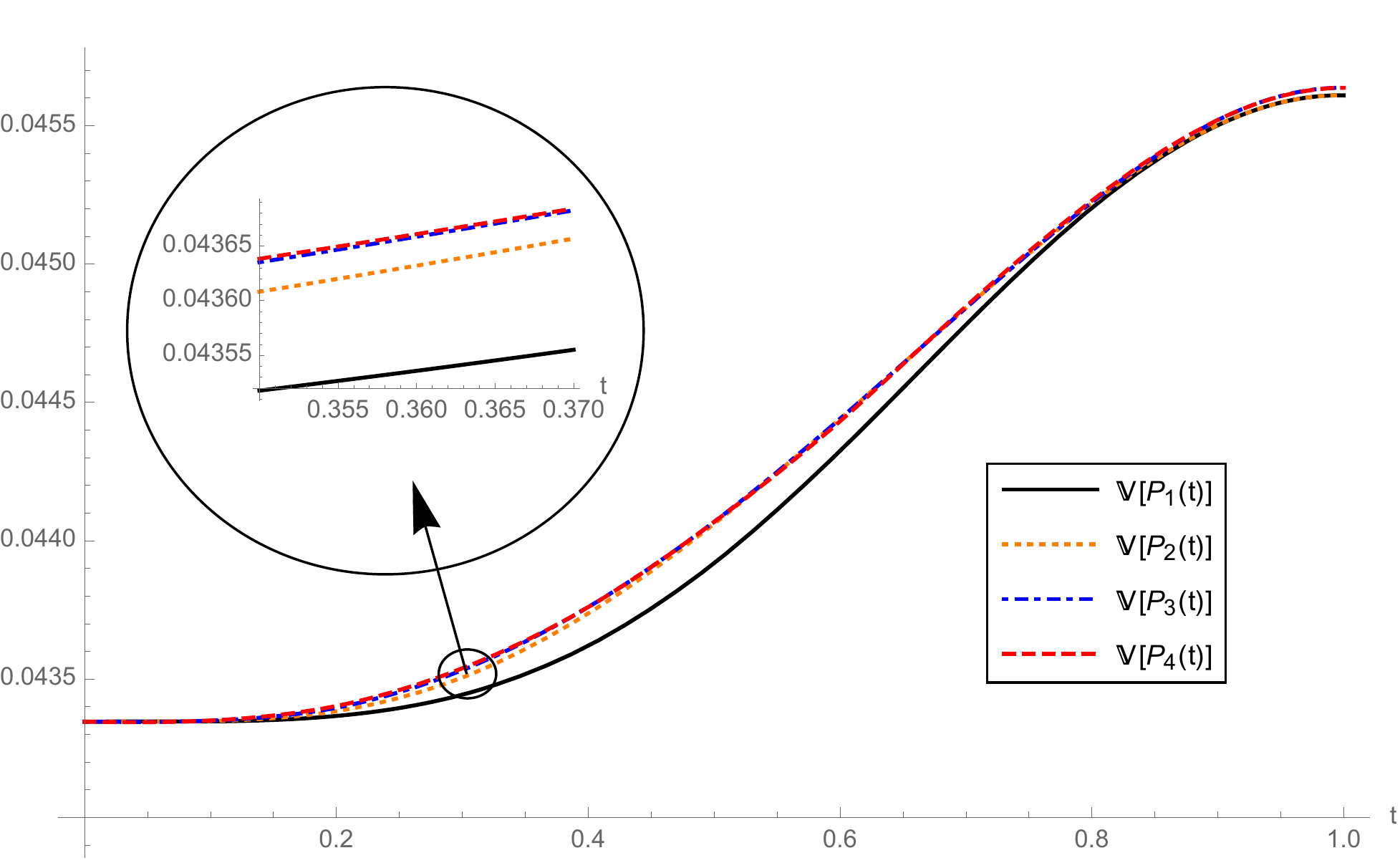}
\end{center}\caption{Example \ref{Ejemplo_Numerico_2}. Left: Approximations of the mean, $\mathbb{E}[P_N(t,\omega)]$. Right: Approximations of the variance, $\mathbb{V}[P_N(t,\omega)]$. In both cases  we have taken  the following orders of truncations $N\in\{1,2,3,4\}$.}\label{grafica-Mean-Variance_ejemplo2}
\end{figure}

\begin{table}[htp]
\begin{center}
\begin{tabular}{*{4}{|c}|}
\hline
$\text{Error}$  & $N=2$ & $N=3$ & $N=4$   \\
\hline
Mean $\hat{e}_N^{\mathbb{E}}$ &  0.000027 & 0.000005 & 0.000002  \\
\hline 
Variance $\hat{e}_N^{\mathbb{V}}$ &   0.000053 & 0.000011 &  0.000004 \\
\hline 
\end{tabular}
\end{center}
\caption{Values of errors $\hat{e}_N^{\mathbb{E}}$ and $\hat{e}_N^{\mathbb{V}}$ for the mean and variance, respectively, given by \eqref{errorMEDIA_VARIANZA_trunca_ej_2} using different orders of  truncation $N\in\{2,3,4\}$, in the context of Example \ref{Ejemplo_Numerico_2}.}
\label{tabla_errores_moments}
\end{table}
\end{example}

\begin{example}\label{Ejemplo_Numerico_3}
We complete the numerical experiments  considering  the random IVP \eqref{logistic_problem_random} on the time interval $\mathcal{T}=[t_0,T]=[-a,a]$, with $a=0.5$. We assume that the initial condition $P_0(\omega)$ has  a Beta distribution truncated to the interval $[0.1,0.9]$ and  parameters $\alpha=7$ and $\beta=10$, i.e. $P_0(\omega) \sim \text{Be}_{[0.1,0.9]} (7;10)$.  Regarding the diffusion coefficient, $A(t,\omega)$,  and in order to apply our theoretical results, we only need to fix the information involved in its KLE \eqref{KLEX}, i.e., a family of zero-mean, unit variance and pairwise  uncorrelated RVs,  $\xi_j(\omega)$, the   mean function, $\mu_A(t)$, and the covariance function, $c_A(s,t)$. Now we will choose:
\begin{itemize}
\item $\xi_j(\omega)$  independent and identically distributed  uniform RVs, $\xi_j(\omega) \sim \text{Un}(-\sqrt{3},\sqrt{3})$. Thus, $\mathbb{E}[\xi_j(\omega) ]=0$, $\mathbb{V}[\xi_j(\omega)]=1$ and $\mathbb{E}[\xi_j(\omega)  \xi_k(\omega)]=0$, if $j \neq k$.
\item Mean function: $\mu_A(t)=0$.
\item Covariance function 
\[
c_A(s,t)=\e^{-c\vert s-t\vert}, \quad (s,t) \in \mathcal{T}\times \mathcal{T},
\]
 where $c=1/b>0$, being $b>0$ the so-called correlation length.
\end{itemize}
According to \cite[pp.~26--29]{spanos}, the eigenvalues $\nu_j$ and eigenfunctions $\phi_j$ of the covariance function $c_A(s,t)$ are given by
\begin{equation}\label{eigen_Ejemplo3}
\begin{array}{ll}
\displaystyle \nu_j=\frac{2c}{w_j^2+c^2},&\displaystyle \phi_j(t)=\frac{\cos (w_j t)}{\sqrt{a+\frac{\sin(2w_j a)}{2w_j}}}, \quad j \text{ odd},\\
\\
\displaystyle \nu_j^*=\frac{2c}{\left(w_j^*\right)^2+c^2},&\displaystyle  \phi_j^*(t)=\frac{\sin (w_j^* t)}{\sqrt{a-\frac{\sin(2w_j^* a)}{2w_j^*}}}, \quad j \text{ even},
\end{array}
\end{equation}
where $w_j$, $w_j^*$  are the solutions of the following transcendental equations
\[
\begin{array}{l}
c-w_j\tan (w_j a )=0, \quad j \text{ odd},\\ \\
w_j^*+c\tan (w_j^* a)=0, \quad j \text{ even}.
\end{array}
\]
Therefore the diffusion SP, $A(t,\omega)$, is represented by the following KLE 
\begin{equation}\label{H4_ejemplo_3_0}
A(t,\omega)=\sum_{j=1}^{\infty} \left( \sqrt{\nu_{2j-1}} \phi_{2j-1} (t) \xi_{2j-1}(\omega)+\sqrt{\nu_{2j}^*} \phi_{2j}^* (t) \xi_{2j}^*(\omega)  \right).
\end{equation}
In order to guarantee that hypothesis \textbf{H2} fulfils, as in the two previous examples we will choose the initial condition $P_0(\omega)$ so that is independent of  the random vector $\boldsymbol{\xi}_{N}(\omega)=(\xi_1(\omega),\ldots,\xi_N(\omega))$, for $N$ arbitrary, but fixed. Now we check that hypotheses  \textbf{H1}, \textbf{H3} and \textbf{H4} hold. First part of hypothesis \textbf{H1} is evident while the second part follows because $\int_{-a}^{a} \mathbb{E}[(A(t,\omega))^2]\, \mathrm{d} t =\int_{-a}^{a}  c_A(t,t)\, \mathrm{d} t=\int_{-a}^{a}  1\, \mathrm{d} t=2a<\infty $. Hypothesis \textbf{H3} can be checked in a similar way as in Example~\ref{Ejemplo_Numerico_1}. To verify that the hypothesis \textbf{H4} holds, we will follow a similar reasoning to the one exhibited in \cite[Remark~2]{KL1}, but now taking advantage of Prop.~\ref{desig_clave}. First, let us observe that using \eqref{K_N} and the independence of $\xi_1(\omega), \ldots,\xi_N(\omega)$, for each $t\in [-a,a]$ and Prop.~\ref{propo_independencia},  one gets,
\begin{equation}\label{H4_ejemplo_3_1}
\begin{array}{ccl}
\mathbb{E}[\e^{2K_N(t,\boldsymbol{\xi}(\omega))}]
&= & \displaystyle
\mathbb{E}
\left[
\exp\left(
2 
\int_{-a}^t
\left(
\sum_{j=1}^N
\left(
\sqrt{\nu_{2j-1}}\,
\phi_{2j-1}(s)\,
\xi_{2j-1}(\omega)
+
\sqrt{\nu_{2j}^{\ast}}\,
\phi_{2j}^{\ast}(s)\,
\xi_{2j}^{\ast}(\omega)
\right)
\right)
\mathrm{d}s
\right)
\right]\\
\\
&= & 
\displaystyle \prod_{j=1}^N
\mathbb{E}
\left[
\exp\left(\lambda_{2j-1}(t)\, \xi_{2j-1}(\omega)\right)
\right]
\mathbb{E}
\left[
\exp \left(\lambda_{2j}^{\ast}(t)\, \xi_{2j}^{\ast}(\omega)\right)
\right],
\end{array}
\end{equation} 
where
\[
\lambda_{2j-1}(t)=2\sqrt{\nu_{2j-1}} \int_{-a}^t \phi_{2j-1}(s)\, \mathrm{d}s,\quad
\lambda_{2j}^{\ast}(t)=2\sqrt{\nu_{2j}^{\ast}} \int_{-a}^t \phi_{2j}^{\ast}(s)\, \mathrm{d}s.
\]
Now, we apply  Prop.~\ref{desig_clave} to each one of the expectations that appear in the last product in \eqref{H4_ejemplo_3_1} (for the first factor we take $\lambda = \lambda_{2j-1}(t)$ and $\xi(\omega)=\xi_{2j-1}(\omega)$,  for the second one, $\lambda = \lambda_{2j}^{\ast}(t)$ and $\xi(\omega)=\xi_{2j}^{\ast}(\omega)$, and $\alpha=-\sqrt{3}$ and $\beta=\sqrt{3}$). This yields
\begin{equation}\label{H4_ejemplo_3_3}
\begin{array}{ccl}
\mathbb{E}[\e^{2K_N(t;\boldsymbol{\xi}(\omega))}]
&\leq & 
\displaystyle
\prod_{j=1}^N
\exp\left(6 \nu_{2j-1} \left( \int_{-a}^t \phi_{2j-1}(s)\, \mathrm{d}s  \right)^2\right)
\exp\left(6 \nu_{2j}^{\ast} \left( \int_{-a}^t \phi_{2j}^{\ast}(s)\, \mathrm{d}s  \right)^2\right)\\
\\
& = &
\exp \left(
\displaystyle 6 \sum_{j=1}^N 
\left(
\nu_{2j-1} 
\left( 
\int_{-a}^t \phi_{2j-1}(s)\, \mathrm{d}s  
\right)^2
+
\nu_{2j} ^{\ast}
\left( 
\int_{-a}^t \phi_{2j}^{\ast}(s)\, \mathrm{d}s  
\right)^2
\right)
\right)
\end{array}
\end{equation} 
Now, we will apply the Cauchy-Schwarz inequality for integrals and the fact that $t\in [-a,a]$, then
\[
\left( \int_{-a}^t \phi_{2j-1}(s)\, \mathrm{d}s  \right)^2
\leq
(t+a) 
 \int_{-a}^t \left( \phi_{2j-1}(s)\right)^2  \mathrm{d}s  
 \leq
2a 
 \int_{-a}^a \left( \phi_{2j-1}(s)\right)^2  \mathrm{d}s , 
\]
and analogously,
\[
\left( \int_{-a}^t \phi_{2j}^{\ast}(s)\, \mathrm{d}s  \right)^2
\leq
2a 
 \int_{-a}^a \left( \phi_{2j}^{\ast}(s)\right)^2  \mathrm{d}s . 
\]
As a consequence, the inequality \eqref{H4_ejemplo_3_3} becomes
\begin{equation}\label{H4_ejemplo_3_4}
\begin{array}{ccl}
\mathbb{E}[
\e^{2K_N(t;\boldsymbol{\xi}(\omega))}]
&\leq &
\exp \left(\displaystyle 
12a
 \sum_{j=1}^{N}
 \left(
\nu_{2j-1} 
 \int_{-a}^a \left( \phi_{2j-1}(s)\right)^2 \mathrm{d}s  
+
\nu_{2j}^{\ast}
 \int_{-a}^a \left( \phi_{2j}^{\ast}(s)\right)^2 \mathrm{d}s 
 \right) 
\right)\\
\\
&= &
\exp\left(\displaystyle 
12a
\left(
\int_{-a}^a 
 \sum_{j=1}^{N}
\nu_{2j-1} 
 \left( \phi_{2j-1}(s)\right)^2 
  \mathrm{d}s 
+
\int_{-a}^a
 \sum_{j=1}^{N}
 \nu_{2j}^{\ast} 
  \left( \phi_{2j}^{\ast}(s)\right)^2 \mathrm{d}s 
 \right) 
\right)\\
\\
&\leq &
\exp\left( \displaystyle 
12a
 \int_{-a}^a 
  \left(
  \sum_{j=1}^{\infty} \nu_{2j-1}  \left( \phi_{2j-1}(s)\right)^2 +  \nu_{2j}^{\ast}  \left( \phi_{2j}^{\ast}(s)\right)^2  \right) 
  \mathrm{d}s 
\right),
\end{array}
\end{equation}
where in the last step we have applied that $\nu_{2j-1}>0$ and $\nu_{2j}^{\ast}>0$ for every $j\geq 1$ (see \eqref{eigen_Ejemplo3}). 

On the other hand,  if we  square the expression of $A(t;\omega)$ given in \eqref{H4_ejemplo_3_0} and afterwards we take the expectation operator and use that $\mathbb{E}[\xi_{2j-1}(\omega)]=\mathbb{E}[\xi_{2j}^{\ast}(\omega)]=0$, $\mathbb{V}[\xi_{2j-1}(\omega)]=\mathbb{V}[\xi_{2j}^{\ast}(\omega)]=1$ and $\mathbb{E}[\xi_j(\omega) \xi_{k}(\omega)]=0$ if $j\neq k$, one obtains
\[
\mathbb{E}[(A(t;\omega))^2]= 
\sum_{j=1}^{\infty} 
\left(
\nu_{2j-1} \left( \phi_{2j-1}(t) \right)^2 
+
\nu_{2j}^{\ast} \left( \phi_{2j}^{\ast}(t) \right)^2
\right) .
\]
Integrating both sides and taking into account that  $A(t,\omega)\in \mathrm{L}^2(\Omega , \mathrm{L}^2(\mathcal{T}))$,  $\mathcal{T}=[-a,a]$ and using the norm defined in \eqref{norma_L_2}, one gets
\[
\int_{-a}^{a} 
\sum_{j=1}^{\infty} 
\left(
\nu_{2j-1} \left( \phi_{2j-1}(t) \right)^2 
+
\nu_{2j}^{\ast} \left( \phi_{2j}^{\ast}(t) \right)^2
\right)
\,\mathrm{d}t 
=
\int_{-a}^{a} 
\mathbb{E}[(A(t;\omega))^2]
\,\mathrm{d}t
=
\left(
\left\| A(t,\omega) \right\|_{\mathrm{L}^2(\Omega , \mathrm{L}^2([-a,a]))}
\right)^2
<\infty
.
\] 
As a consequence, using this last conclusion in expression \eqref{H4_ejemplo_3_4}, one derives that $\mathbb{E}[
\e^{2K_N(t;\boldsymbol{\xi}(\omega))}]<\infty$ for every $t\in [-a,a]$ and for all $N\geq 1$ positive integer.  Therefore, the hypothesis \textbf{H4} fulfils.

Now, we do not show the explicit algebraic expression of the approximations, $f_1^N(p,t)$, because it is somewhat cumbersome, but is clear that it could be calculated in the same manner  we did in the two previous examples. In Figure \ref{grafica1PDFs_ejemplo3}, we   show the surfaces corresponding to those  approximations  for $N=1$ and $N=2$. From these two plots, we can observe that both surfaces are very similar, then showing a fast convergence. In Figure \ref{grafica1PDFs-t_ejemplo3}, we show the approximations $f_1^N(p,t)$ at different time instants $t\in\{-0.25, 0, 0.25  \}$ and using different orders of truncations $N\in\{1,2,3\}$. Again, we observe  fast convergence in all these cases. We use an analogous  measure to the one  defined in \eqref{error-PDF_trunca}, i.e.,
\begin{equation}\label{error-PDF_trunca_ej_3}
\hat{e}_N^{PDF}(t)=\int_{-0.5}^{0.5} \left|  f_1^N(p,t)-f_1^{N-1}(p,t) \right| \mathrm{d}p, \quad N=2,3,\ldots,
\end{equation}
to illustrate this convergence. In Table \ref{tabla_errores_ejemplo_3}, we have collected figures of $\hat{e}_N^{\text{PDF}}(t)$ for the values of $t$ and $N$ previously indicated.

\begin{figure}[htp]
\begin{center}
\includegraphics[width=0.49\textwidth]{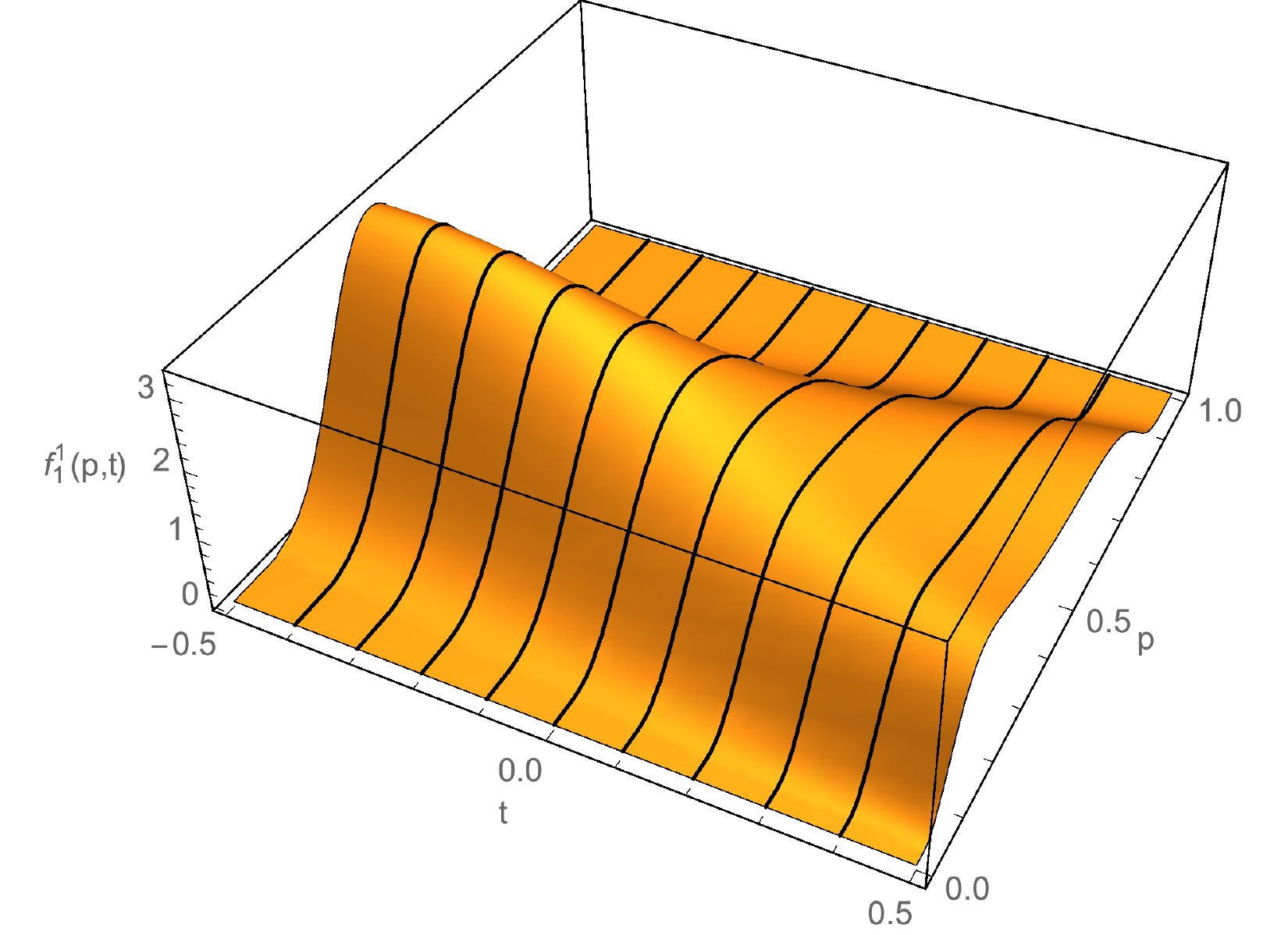}
\includegraphics[width=0.49\textwidth]{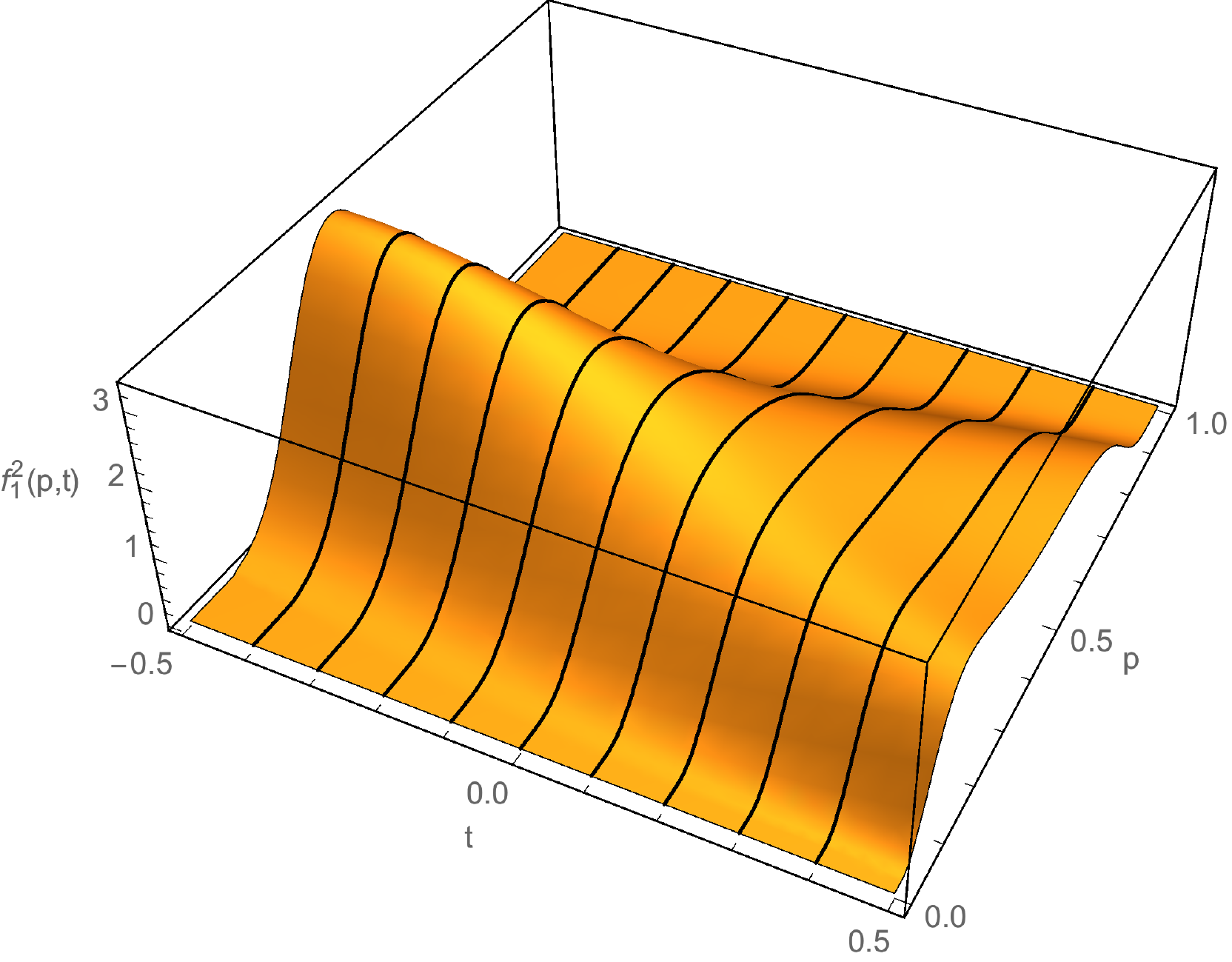}
\end{center}\caption{Example  \ref{Ejemplo_Numerico_3}. Surfaces of the 1-PDF, $f_1^N(p,t)$,  for  $N=1$ (Left) and $N=2$ (Right).}\label{grafica1PDFs_ejemplo3}
\end{figure}

\begin{figure}
\begin{center}
\includegraphics[width=0.32\textwidth]{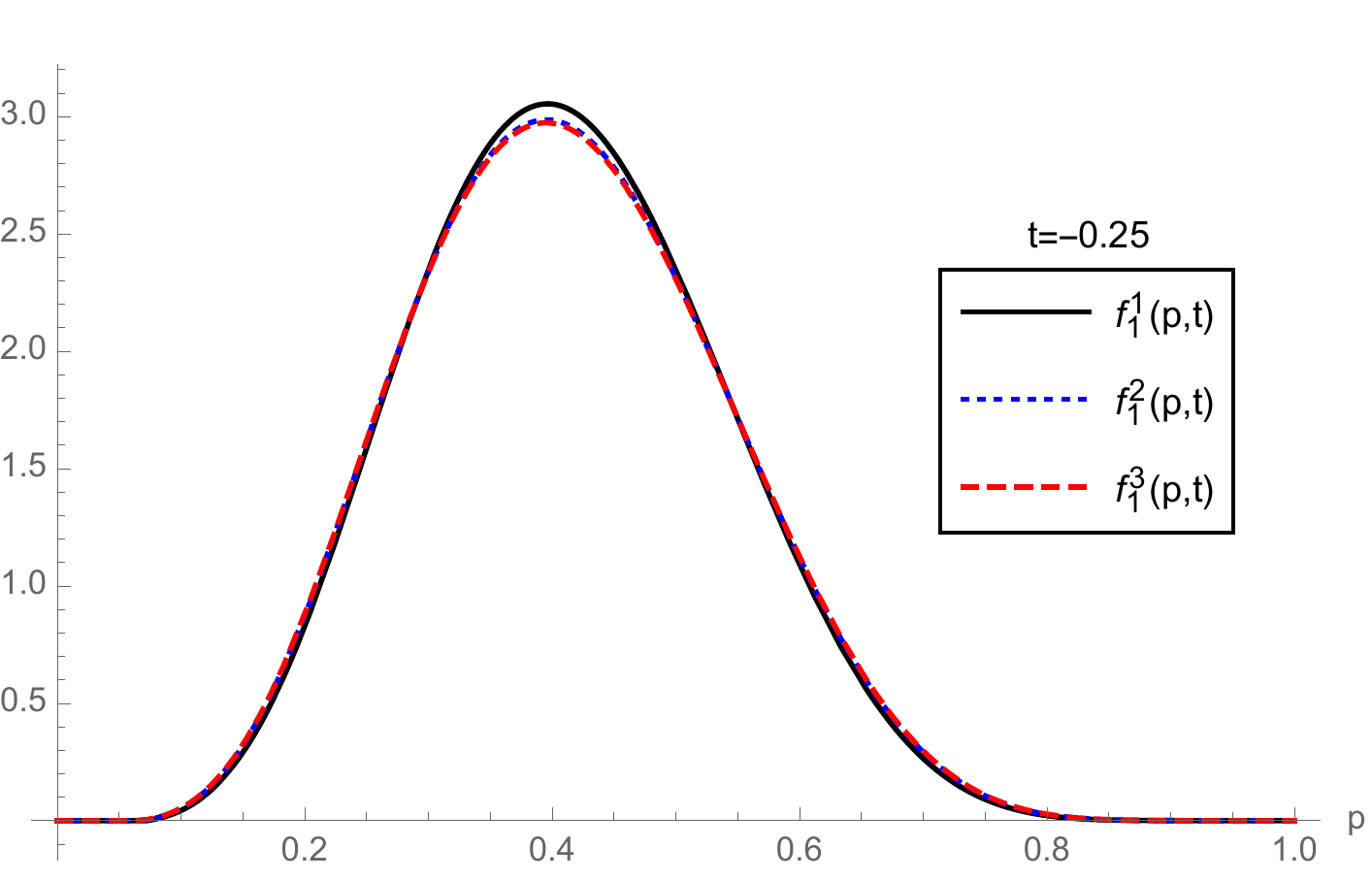}
\includegraphics[width=0.32\textwidth]{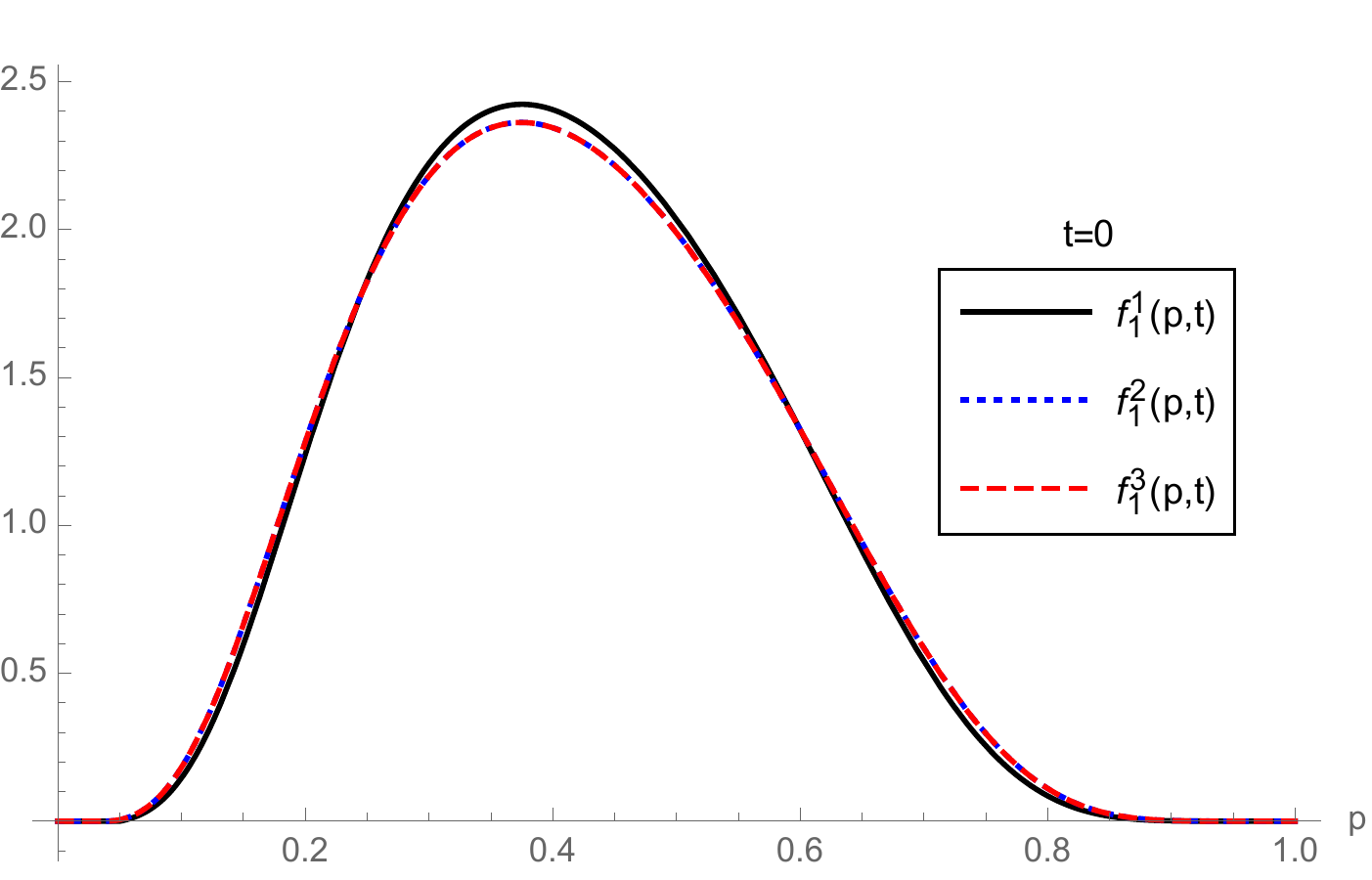}
\includegraphics[width=0.32\textwidth]{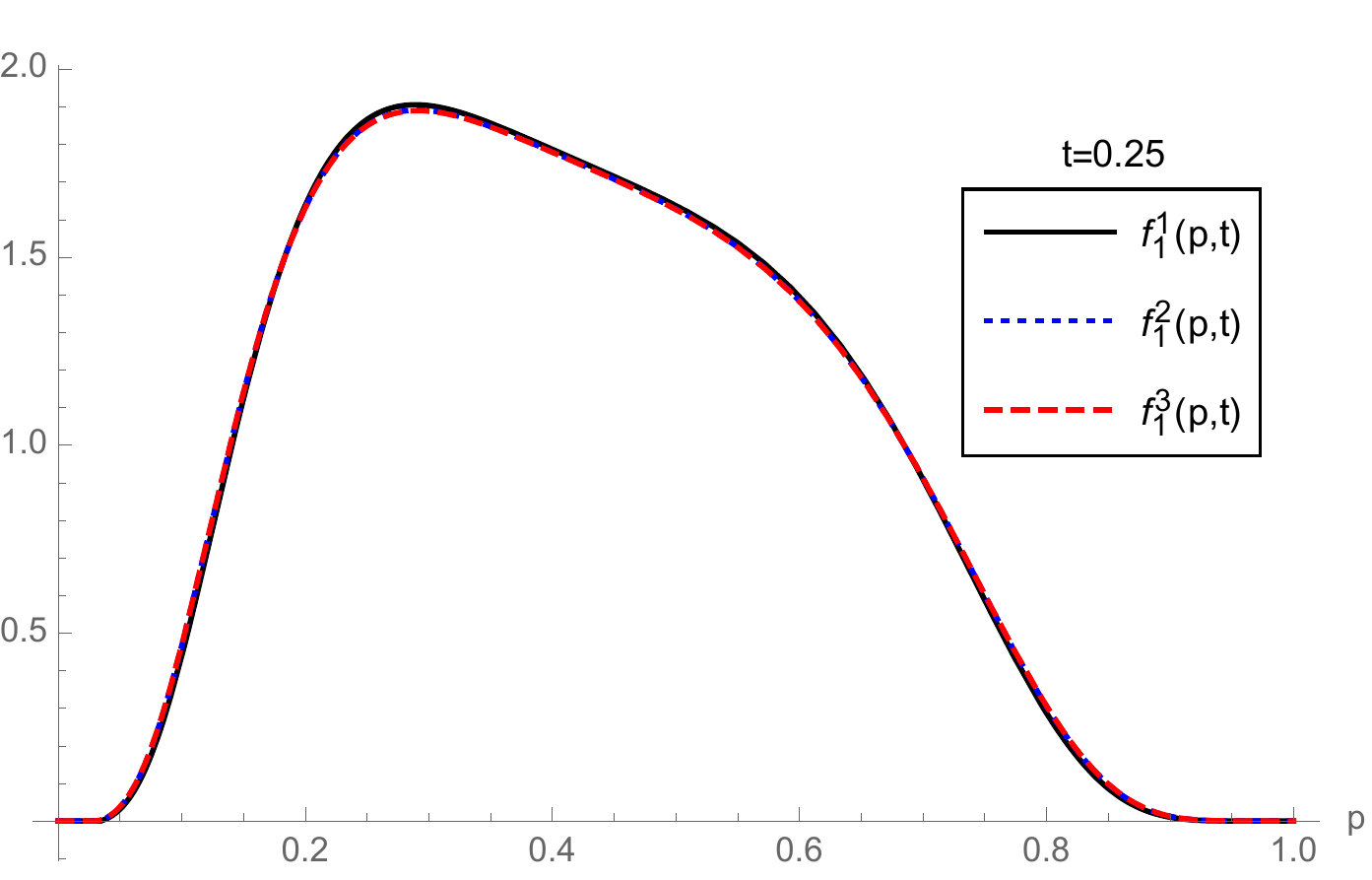}
\end{center}\caption{Example \ref{Ejemplo_Numerico_3}. Curves of the 1-PDF, $f_1^N(p,t)$ at  three different times $t=-0.25$ (Left), $t=0$ (Center) and $t=0.25$ (Right) using, in each case, different orders of truncations $N\in \{ 1,2,3\}$.}\label{grafica1PDFs-t_ejemplo3}
\label{grafica1PDFs-t_ejemplo3}
\end{figure}

\begin{table}[htp]
\begin{center}
\begin{tabular}{|l|c|c|}
\hline
$\hat{e}_N^{\text{PDF}}(t)$  & $N=2$ & $N=3$    \\
\hline
$t=-0.25$ &    0.022077 & 0.004105 \\
\hline 
$t=0$ &   0.029739 & 0.000044  \\
\hline 
$t=0.25$ &  0.009479 & 0.000975\\
\hline 
\end{tabular}
\end{center}
\caption{Error measure $\hat{e}_N^{\text{PDF}}(t)$ defined by \eqref{error-PDF_trunca_ej_3} for different time instants, $t\in \{ -0.25,0,0.25 \}$, and truncation orders $N\in\{2,3\} $, in the context of Example 3.}
\label{tabla_errores_ejemplo_3}
\end{table}

We finally compute approximations of the mean, $\mathbb{E}[P_N(t,\omega)]$,  and the variance, $\mathbb{V}[P_N(t,\omega)]$, for  $N\in\{1,2,3,4\}$. These approximations have been represented in Figure \ref{grafica-Mean-Variance_ejemplo3}. From these graphical representations we  evince fast convergence of both statistical moments.  In Table \ref{tabla_errores_moments_ejemplo_3}, we illustrate numerically this convergence by computing  the total difference  between consecutive approximations of the mean and  variance on the whole time interval via the  expressions given in \eqref{errorMEDIA_VARIANZA_trunca_ej_3}.
\begin{equation}\label{errorMEDIA_VARIANZA_trunca_ej_3}
\hat{e}_N^{\mathbb{E}}=\int_{-0.5}^{0.5} \left| \mathbb{E}[P_{N}(t,\omega)]-\mathbb{E}[P_{N-1}(t,\omega)]\right| \mathrm{d} t,\qquad
\hat{e}_N^{\mathbb{V}}=\int_{-0.5}^{0.5} \left| \mathbb{V}[P_{N}(t,\omega)]-\mathbb{V}[P_{N-1}(t,\omega)]\right| \mathrm{d} t.
\end{equation}

\begin{figure}
\begin{center}
\includegraphics[width=0.49\textwidth]{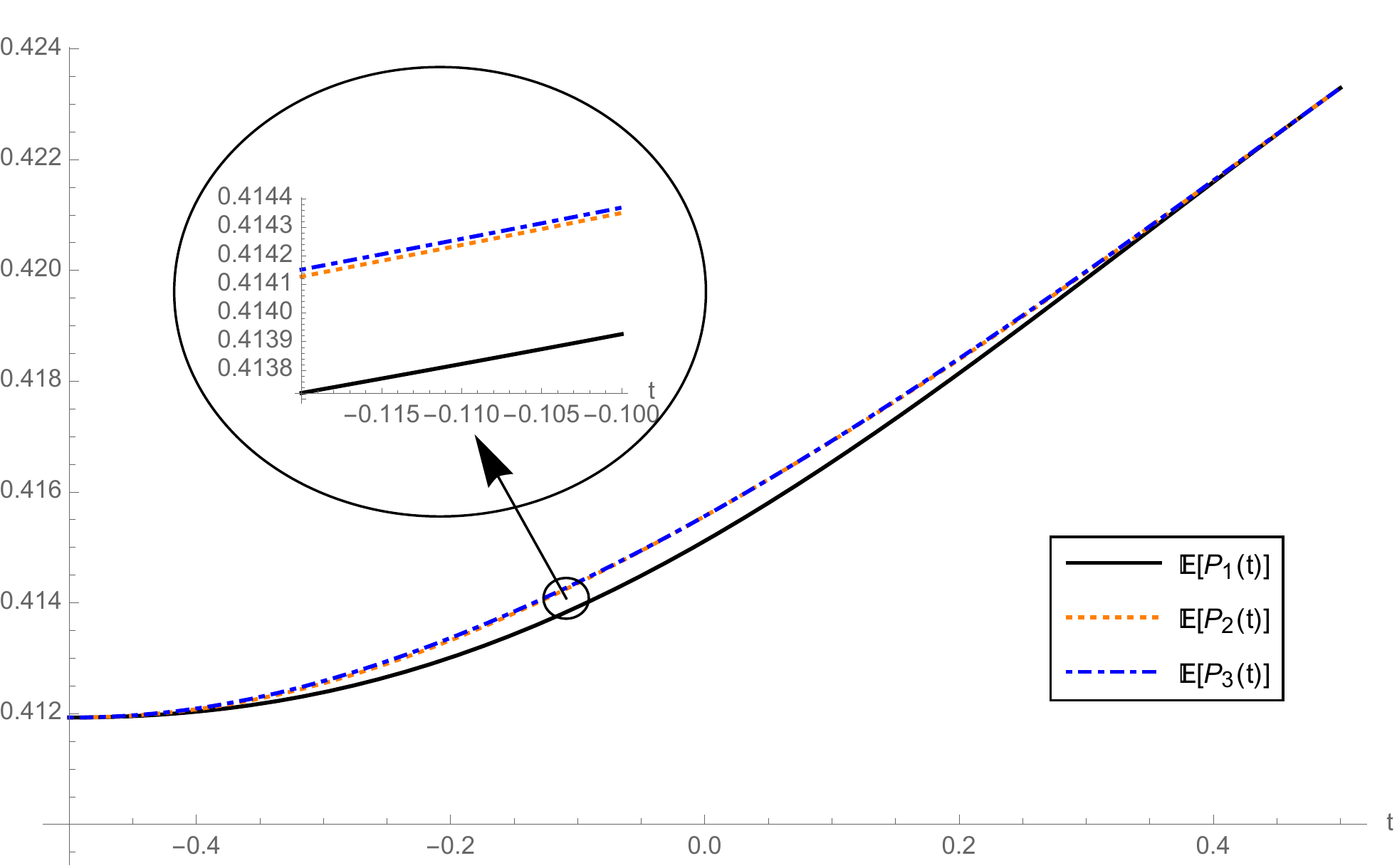}
\includegraphics[width=0.49\textwidth]{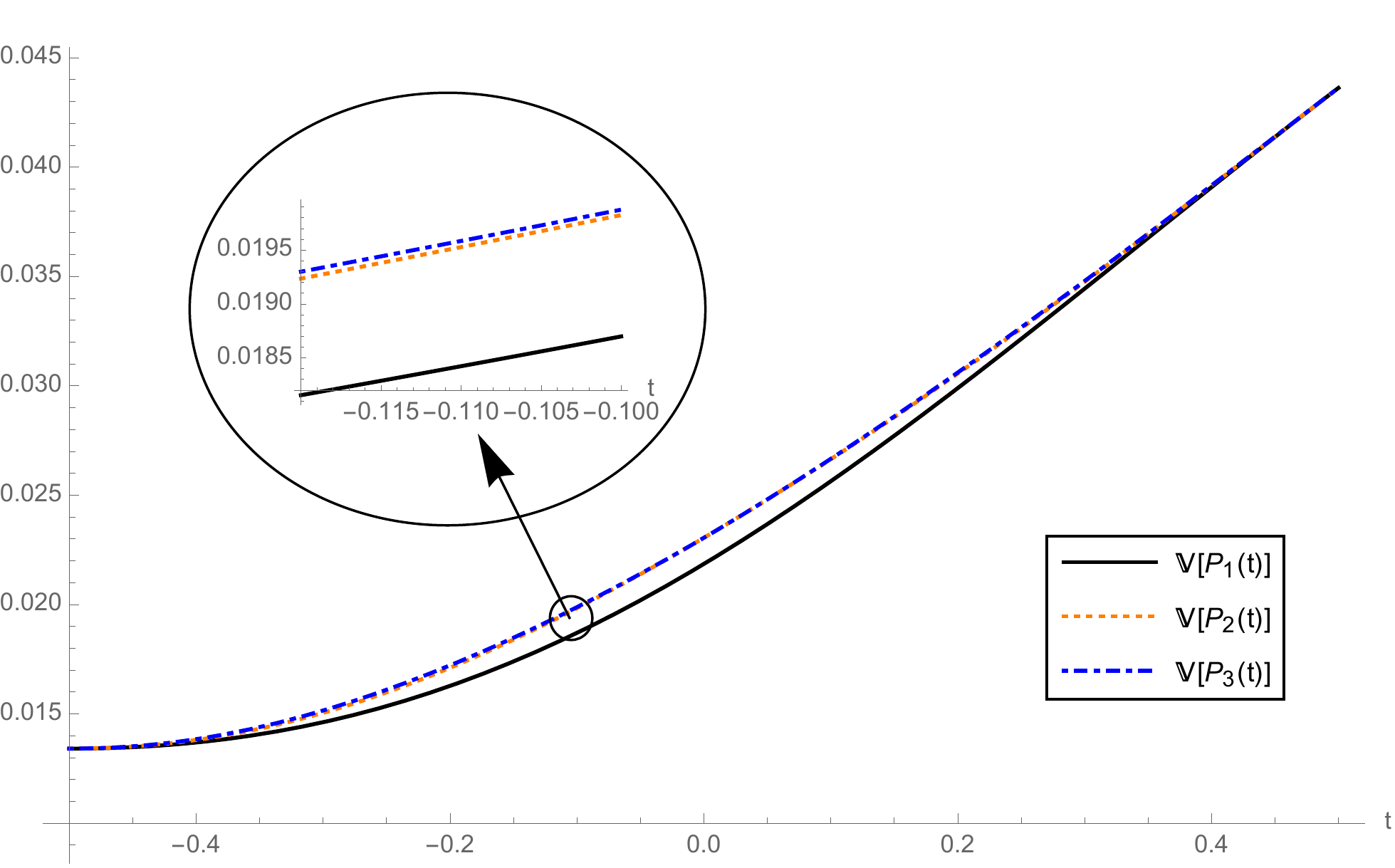}
\end{center}\caption{Example \ref{Ejemplo_Numerico_3}. Left: Approximations of the mean, $\mathbb{E}[P_N(t,\omega)]$. Right: Approximations of the variance, $\mathbb{V}[P_N(t,\omega)]$. In both cases  we have taken  the following orders of truncations $N\in\{1,2,3\}$.}\label{grafica-Mean-Variance_ejemplo3}
\end{figure}

\begin{table}[htp]
\begin{center}
\begin{tabular}{*{3}{|c}|}
\hline
$\text{Error}$  & $N=2$ & $N=3$   \\
\hline
Mean $\hat{e}_N^{\mathbb{E}}$ &  0.000216 & 0.000016 \\
\hline 
Variance $\hat{e}_N^{\mathbb{V}}$ &   0.000575 & 0.000042  \\
\hline 
\end{tabular}
\end{center}
\caption{Values of errors $\hat{e}_N^{\mathbb{E}}$ and $\hat{e}_N^{\mathbb{V}}$ for the mean and variance, respectively, given by \eqref{errorMEDIA_VARIANZA_trunca_ej_3} using different orders of  truncation $N\in\{2,3\}$, in the context of Example \ref{Ejemplo_Numerico_3}.}
\label{tabla_errores_moments_ejemplo_3}
\end{table}
\end{example}
\section{Conclusions}\label{conclusiones}
In this paper we have studied a generalization of the random logistic differential equation  consisting of assuming that the diffusion coefficient is a stochastic process and with a random initial condition. Under  general hypotheses on random data, we have constructed approximations of the first probability density function of the solution stochastic process. The key tools for conducting our analysis have been the Random Variable Transformation method together with the Karhunen-Lo\`{e}ve expansions. Our theoretical findings have been illustrated by means of several examples. To the best of our knowledge, it is first time that our approach is applied to a random non-autonomous nonlinear  differential equation. We think that this contribution  can be  useful to study other important random nonlinear differential equations.    

\section*{Acknowledgements}
This work has been partially supported by the Ministerio de Econom\'{i}a y Competitividad grant MTM2017-89664-P. Ana Navarro Quiles acknowledges the postdoctoral contract financed by DyCon project funding from the European Research Council (ERC) under the European Union’s Horizon 2020 research and innovation programme (grant agreement No 694126-DYCON).

The authors express their deepest thanks and respect to the editor and reviewers for their valuable comments.

\section*{Conflict of Interest Statement}
The authors declare that there is no conflict of interests regarding the publication of this article.




\section*{References}

\bibliographystyle{elsarticle-num}
\bibliography{rvt}







\end{document}